\documentclass[journal]{IEEEtran}
\ifCLASSINFOpdf
\else
\fi
\usepackage{amsfonts}
\usepackage{bm}
\usepackage{amsmath}
\usepackage{amssymb}
\usepackage{graphicx} 
\usepackage{epstopdf}
\usepackage{float}
\usepackage{makecell}
\usepackage{algorithmic}
\usepackage{multirow}
\usepackage{cite}
\usepackage{subfigure}
\usepackage{epsfig} 
\usepackage{color}
\usepackage[ruled,vlined]{algorithm2e}

\begin{document}

%
\title{Collision-Free Trajectory Design for 2D Persistent Monitoring Using Second-Order Agents}
%
%
%
\author{Yan-Wu Wang,~\IEEEmembership{Senior Member,~IEEE,}
        Ming-Jie Zhao, Wu Yang,
        Nan Zhou,~\IEEEmembership{Student Member, IEEE}
        and\\ Christos G. Cassandras,~\IEEEmembership{Fellow,~IEEE}
\thanks{This work is supported by the National Natural Science Foundation of China under Grants 61773172 and 51537003, the Natural Science Foundation of Hubei Province of China (2017CFA035), the Fundamental Research Funds for the Central Universities (2018KFYYXJJ119) and the Program for HUST Academic Frontier Youth Team.}

\thanks{Y.-W. Wang, M.-J. Zhao and  W. Yang are with the School of Artificial Intelligence and Automation, Huazhong University of Science and Technology, Wuhan, 430074,  China, and also with the Key Laboratory of Image Processing and Intelligent Control (Huazhong University of Science and Technology), Ministry of Education, China (e-mail: wangyw@hust.edu.cn; zhaomj@hust.edu.cn; wuyangac@hust.edu.cn; biwuyang@163.com).}
\thanks{N. Zhou is with the Division of Systems Engineering, Boston University, Boston, MA 02446 USA (e-mail: nanzhou@bu.edu).}
\thanks{C.G. Cassandras is with the Division of Systems Engineering and Department of Electrical and Computer Engineering, Boston University,  Boston, MA 02446 USA (e-mail: cgc@bu.edu).}}
\maketitle

\begin{abstract}
This paper considers a two-dimensional persistent monitoring problem by controlling movements of second-order agents to minimize some uncertainty metric associated with targets in a dynamic environment. In contrast to common sensing models depending only on the distance from a target, we introduce an active sensing model which considers the distance between an agent and a target, as well as the agent's velocity. We propose an objective function which can achieve a collision-free agent trajectory by penalizing all possible collisions. Applying structural properties of the optimal control derived from the Hamiltonian analysis, we limit agent trajectories to a simpler parametric form under a family of 2D curves depending on the problem setting, e.g. ellipses and Fourier trajectories. Our collision-free trajectories are optimized through an event-driven Infinitesimal Perturbation Analysis (IPA) and gradient descent method. Although the solution is generally locally optimal, this method is computationally efficient and offers an alternative to other traditional time-driven methods. Finally, simulation examples are provided to demonstrate our proposed results.
\end{abstract}


%
\IEEEpeerreviewmaketitle

\section{Introduction}\label{ZMJ_PersistentMonitoring_section_1}
%
%
%
%
\IEEEPARstart{I}{n} recent years, attention has been drawn to persistent monitoring because of its wide range of applications such as ocean sampling \cite{paley2008cooperative,smith2011persistent}, city surveillance \cite{michael2011persistent,da2017unmanned}, and traffic monitoring \cite{puri2005survey,vian2014traffic}. Persistent monitoring often arises in a large dynamically changing environment which requires agents to move to cover the monitoring region. The fundamental problem of persistent monitoring is to design suitable motion strategies for agents to meet the monitoring requirements cooperatively. Up to now, there have been many studies on persistent monitoring in both 1D \cite{zhou2016optimal}, 2D \cite{tokekar2015visibility} and 3D \cite{yu2016correlated} spaces. For the 1D problems, \cite{zhou2018optimal} recently addresses a 1D persistent monitoring problem which only involves a finite number of targets and designed optimal trajectories for agents. \cite{wang2018optimal1} introduces a monitoring index to differentiate regions with different monitoring importance. Moreover, \cite{wang2018optimal2} solves the 1D persistent monitoring utilizing second-order agents with physical constraints, e.g. bounded acceleration and velocity.

Compared with the 1D problems, the 2D monitoring problems find more practical applications; the control strategies are more complicated; and the optimal solutions are more challenging to obtain. \cite{lin2015optimal} restricts agent trajectories to ellipses and proved that they outperform the linear ones under certain conditions in a 2D mission space. \cite{song2014optimal} studies monitoring finite stationary targets distributed nearby a given closed trajectory in 2D space. And \cite{smith2012persistent} presents optimal speed controllers along a given closed path to stabilize a changing environment.

Different from 1D and 2D problems, the key to solving 3D monitoring is not to focus on trajectories of agents, but to model it as a scheduling or visiting problem. \cite{stump2011multi} casts a persistent surveillance problem as a Vehicle Routing Problem with Time Windows, and develops a locally optimal path within a set time horizon, which then is executed repeatedly. \cite{yu2015persistent} designs an optimal scheduling scheme for one robot tasked to monitor several events that are occurring at different locations.

The aforementioned literature deals with the persistent monitoring problem without obstacles in the mission space. However, in many applications, the existence of obstacles restricts the movements of agents, thus exacerbating the complexity of the problem. \cite{wang2018optimal2} proposes a collision avoidance algorithm by repeatedly adjusting the designed agent trajectories. And \cite{soltero2011collision} solves collision problems in 2D monitoring through stopping policies on given closed paths. Further, in related research on coverage control \cite{song2013persistent}, there are some results involving avoiding obstacles. Specifically, \cite{franco2015persistent} achieves collision avoidance through a bounded repulsive avoidance control law. \cite{palacios2017optimal} keeps robots at a safety distance from obstacles through speed functions. However, such methods of avoiding obstacles are not applicable in persistent monitoring. Thus the design of collision-free trajectories is still an open question.

Motivated by the above discussion, we will first introduce the persistent monitoring problem using second-order agents and then consider all possible collisions (agent to agent, agent to obstacle). Our goal is to $i)$ find a solution to an optimal control problem for second-order agents visiting targets with different weights, while also avoiding all possible types of collisions. Through standard optimal control theoretic techniques \cite{bryson1975applied}, we show that agents move with maximal accelerations; $ii)$ translate the optimal control problem into a tractable parametric optimization problem. Inspired by \cite{lin2015optimal}, we represent agent trajectories using various families of parametric curves. $iii)$ determine optimal parametric agent trajectories under the parametric families used. IPA \cite{cassandras2010perturbation} is utilized for the gradient evaluation and the gradients are then used for optimization through gradient descent methods. Moreover, the IPA gradient is updated in an event-driven manner, hence not requiring continuous time-driven updates. Note that our algorithm optimizes the agent trajectories while ensuring agent safety by penalizing all possible collisions.

The main contributions of this paper can be summarized as follows:
\begin{itemize}
\item Second-order agents are utilized to execute 2D monitoring task, which brings our treatment of persistent monitoring one step closer to realistic applications. Moreover, bounds with respect to agent accelerations and velocities are both considered.
\item Targets with different weight coefficients are considered, which is different from \cite{lin2015optimal} where all targets are homogeneous.
\item A new sensing model is constructed, which is quite different from those limited to dependence only on the distance between an agent and a target in \cite{zhou2016optimal,zhou2018optimal,wang2018optimal1,wang2018optimal2,lin2015optimal}. Our sensing model depends not only on relative distances between agents and targets but also on agents' velocities \cite{adlakha2003critical,zang2006critical}, and finds use in applications such as vision-based monitoring \cite{baseggio2010distributed,carli2011distributed}.
\item Obstacles in the mission space are considered. Thus, we propose a novel objective function which leads to a collision-free solution by penalizing all possible collisions.
\end{itemize}

The rest of the paper is organized as follows. In Section \ref{ZMJ_PersistentMonitoring_section_2}, an optimal control formulation for the persistent monitoring problem is proposed. In Section \ref{ZMJ_PersistentMonitoring_section_3},  some preliminary insights about the solution are given through Hamiltonian analysis. Section \ref{ZMJ_PersistentMonitoring_section_4} parameterizes the agent trajectories under a selected family of 2D parametric curves and optimize the parameters using an IPA-based gradient descent algorithm. Section \ref{ZMJ_PersistentMonitoring_section_5} provides simulations and Section \ref{ZMJ_PersistentMonitoring_section_6} concludes the paper.

\section{Persistent Monitoring Problem Formulation}\label{ZMJ_PersistentMonitoring_section_2}
\begin{figure}
\centering
\includegraphics[height=3cm,width=6cm]{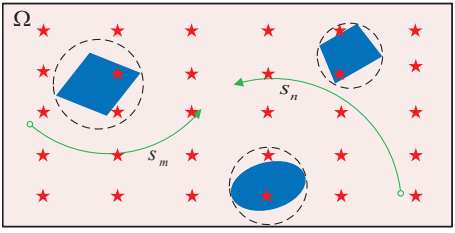}
\caption{Persistent monitoring mission space. Blue areas are obstacles. Red stars represent target points. The curves $s_{m}$ and $s_{n}$ represent feasible agent trajectories.} \label{ZMJ_PersistentMonitoring_Environment}
\end{figure}
Consider a 2D rectangular mission space $\bm{\Omega} \equiv [0, L_{1}] \times [0, L_{2}] \subset \mathbb{R}^{2}$ as shown in Fig. \ref{ZMJ_PersistentMonitoring_Environment}. Select $\{z_{i} = [z_{i}^{x}, z_{i}^{y}]  \in \bm{\Omega}, i=1,...,M\}$ as target points to be monitored. Note that some target points may be located within obstacles. We assume sensors can monitor these points without going into obstacles. This can model, for example, sensors located outside buildings which can detect targets inside the buildings.
\subsection{Agent Dynamics}
We consider $N$ cooperating agents assigned to monitor stationary target points over the time horizon $[0, T]$. The dynamics of agent $n$ can be described by
\begin{equation}
\begin{aligned}
\dot{x}_{n}(t)=f(x_{n},U_{n})=\begin{bmatrix} \bm{{\rm 0}}_{2,2} & \bm{{\rm I}}_{2} \\ \bm{{\rm 0}}_{2,2} & \bm{{\rm 0}}_{2,2} \end{bmatrix}x_{n}(t)
+ \begin{bmatrix} \bm{{\rm 0}}_{2,2} \\ \bm{{\rm I}}_{2} \end{bmatrix}U_{n}(t)
\end{aligned}
\end{equation}
where $n = 1,...,N$, $\bm{{\rm 0}}_{2,2}$ is a $2 \times 2$ zero matrix, $\bm{{\rm I}}_{2}$ is a $2 \times 2$ identity matrix and $U_{n}$ is the control input.
The state $x_{n}(t) = [s_{n}(t), v_{n}(t)]^{\top}$, where $s_{n}(t)=[s_{n}^{x}(t), s_{n}^{y}(t)]^{\top}$ and $v_{n}(t)=[v_{n}^{x}(t), v_{n}^{y}(t)]^{\top}$ represents the position  and the velocity of the agent $n$ at time $t$ respectively. In addition, we define the control input $U_{n}(t) = [u_{n}(t)\cos\theta_{n}(t), u_{n}(t)\sin\theta_{n}(t)]^{\top}$, where $u_{n}(t)$ is the magnitude of the acceleration and $\theta_{n}(t)$ is the agent's heading that satisfies $0 \leq \theta_{n}(t) < 2\pi$. Without loss of generality, we assume that the velocity and the acceleration of each agent are bounded such that
\begin{equation}\label{ZMJ_PersistentMonitoring_v_constraint}
\begin{split}
0 \leq \| v_{n}(t)\| \leq  v_{n}^{max}, \; n=1,...,N
\\
0 \leq \Vert u_{n}(t) \Vert \leq u_{n}^{max}, \; n = 1,...,N
\end{split}
\end{equation}
where $v_{n}^{max}$ is subject not only to the maximum power constraint \cite{wang2018optimal2}, but also to the performance of the sensor aboard the agent (which will be mentioned in the following).

\textbf{Remark 1.} Note that the control input $u_{n}(t)$ here is the acceleration of the agent $n$ as opposed to the velocity in \cite{lin2015optimal}. Moreover, as will be seen in the next section, following the centralized determination of optimal trajectories, the monitoring task is executed in an open-loop fashion by each agent.

\subsection{Agent Sensing Model}
According to \cite{adlakha2003critical,zang2006critical}, the detection probability of a sensor is inversely proportional to the movement velocity. Thus, here we assume that the actual sensing strength of a sensor is determined by the distance between a sensor and a target, as well as its velocity. Specifically, a closer distance and a smaller velocity lead to higher sensing strength. The sensing effectiveness is zero when the distance exceeds a finite sensing range $r_{n}$ or the velocity of the agent exceeds a certain threshold $\beta_{n}$ $(> v_{n}^{max}$ to avoid meaningless agent movements). Define a function $p_{n}(z_{i})$ which represents the probability that an event at location $z_{i} = [z_{i}^{x}, z_{i}^{y}]$ is detected by agent $n$. In addition to those properties of the sensing model considered in \cite{lin2015optimal}, $p_{n}(z_{i})$ is also monotonically non-increasing function in the agent's velocity.
\begin{equation}\label{ZMJ_PersistentMonitoring_sensingmodel}
\begin{aligned}
p_{n}(z_{i}) =  \left\{
\begin{array}{rcl}
(1 - \frac{D(z_{i},s_{n})}{r_{n}})(1 - \frac{\|v_{n}\|}{\beta_{n}}),\qquad\quad\;\;\\
{{\rm if}\; D(z_{i},s_{n}) \leq r_{n}\;{\rm and}\;\|v_{n}\| \leq \beta_{n}}\\
\\
0, {{\rm if}\; D(z_{i},s_{n}) > r_{n}\;{\rm or}\;\|v_{n}\| > \beta_{n}}
\end{array} \right.
\end{aligned}
\end{equation}
where $D(z_{i},s_{n}) = \| z_{i} - s_{n} \|$ is the Euclidean distance between the agent $n$ and the target position $z_{i}$, and $v_{n}$ is the agent velocity (we assume that targets are stationary).

Since there may be multiple agents, the joint probability that an event occurring at $z_{i}$ is detected, denoted by $P_{i}(\mathbf{x}(t))$, is given by
\begin{equation}\label{ZMJ_PersistentMonitoring_probability}
\begin{aligned}
P_{i}(\mathbf{x}(t)) = 1 - \prod_{n = 1}^{N}[1 - p_{n}(z_{i})]
\end{aligned}
\end{equation}
where $\mathbf{x}(t) = [\mathbf{s}(t),\mathbf{v}(t)]^{\top}$ with $\mathbf{s}(t) = [s_{1}(t), ..., s_{N}(t)]$ and $\mathbf{v}(t) = [v_{1}(t),...,v_{N}(t)]$.

\subsection{Target Dynamics}
We associate an uncertainty function $R_{i}(t)$ with every target point $z_{i}$, which possesses similar properties to the model in \cite{lin2015optimal}: (\emph{i}) $R_{i}(t)$ increases with a prespecified rate $A_{i}$ if $P_{i}(t) = 0$, i.e. there is no agent detecting target $z_{i}$; (\emph{ii}) $R_{i}(t)$ decreases with a fixed rate $B$ if $P_{i}(t) = 1$; (\emph{iii}) the decrease of $R_{i}(t)$ is proportional to the joint probability $P_{i}(\mathbf{x}(t))$; (\emph{iv}) $R_{i}(t) \geq 0$ for all \emph{t}. Thus, the dynamics of $R_{i}(t)$ are
\begin{equation}\label{ZMJ_PersistentMonitoring_Uncertainty}
\begin{aligned}
\dot{R}_{i}(t) =  \left\{
\begin{split}
& 0, \qquad {{\rm if}\; R_{i}(t) = 0,\; A_{i} \leq BP_{i}(\mathbf{x}(t))}\\
& A_{i} - BP_{i}(\mathbf{x}(t)), \qquad\qquad{{\rm otherwise}}
\end{split} \right.
\end{aligned}
\end{equation}
We further assume that initial conditions $R_{i}(0), i=1,...,M$, are given and that $B > A_{i} > 0$ for stability.

\subsection{Optimal Control Problem}
Our goal is to minimize the uncertainty accumulated across all target points. We define $J_1(t)$ to be the weighted sum of target uncertainties:
\begin{equation}\label{ZMJ_PersistentMonitoring_J1}
\begin{aligned}
J_{1}(t) = \displaystyle{\sum_{i = 1}^M}\sigma_{i}R_{i}(t)
\end{aligned}
\end{equation}
The weight coefficients $\sigma_{i}$ are set to capture the relative importance of different targets. Note that the problem setting in \cite{lin2015optimal}, where $\sigma_{i} = 1, \text{ for } i=1,...,M$, is a special case of our setting here.

Moreover, in contrast to \cite{lin2015optimal} where each agent is represented as a point mass and collisions among agents are ignored, we will consider the sizes of agents in this work. Note that to avoid collisions in persistent monitoring settings, any two agents cannot share the same location at the same time instant. Considering the size of each agent, for agent $n$ we define a safety radius $\rho_{n} > 0$, and the corresponding safety disk $Q_{n} = \left\{ x\in\bm{\Omega} \bm{\mid} \| x-s_{n}(t) \| \leq \rho_{n}\right\}$. We consider that a collision occurs between agents $p$ and $q$ at some location only if $Q_{p} \cap Q_{q} \neq \emptyset$. Obviously, to avoid agent collisions, we must ensure that the Euclidean distance $d_{pq}(t) = \| s_{p}(t) - s_{q}(t) \| \geq \rho_{p} + \rho_{q}$ at all times. To capture the collisions among agents, we define
\begin{equation}\label{ZMJ_PersistentMonitoring_dpq}
\begin{aligned}
d_{pq}^{-}(t) = {\rm min}(0, d_{pq}(t) - \rho_{p} - \rho_{q})
\end{aligned}
\end{equation}
First, $d_{pq}^{-}(t) \leq 0$ for all $t \in [0,T]$. Second, a collision-free trajectory satisfies $d_{pq}^{-}(t) = 0$ for all $t$. Considering all possible collisions among agents, we define the agent collision avoidance component $J_2(t)$ of the objective function as follows:
\begin{equation}\label{ZMJ_PersistentMonitoring_J2}
\begin{aligned}
J_{2}(t) = \displaystyle{\sum_{q = p+1}^N} \displaystyle{\sum_{p = 1}^{N-1}}d_{pq}^{-}(t)
\end{aligned}
\end{equation}

In addition, the presence of obstacles in the mission space (as shown in Fig. \ref{ZMJ_PersistentMonitoring_Environment}) usually restricts the feasible trajectories. This brings another level of difficulty to the optimization problem. For simplicity, we ignore shapes of the obstacles by covering them using circumscribed circles (shown as dashed circles in Fig. \ref{ZMJ_PersistentMonitoring_Environment}) whose centers are denoted by $\omega_{l} = [\omega_{l}^{x}, \omega_{l}^{y}]$ and radii by $r_{l}, l = 1,...,L$.
Conservatively, we can ensure the safety of agent $n$ by letting the Euclidean distance $d_{ln}(t) = \| \omega_{l} - s_{n}(t) \| \geq r_{l} + \rho_{n}$. Similar to \eqref{ZMJ_PersistentMonitoring_dpq}, we define
\begin{equation}\label{ZMJ_PersistentMonitoring_dln}
\begin{aligned}
d_{ln}^{-}(t) = {\rm min}(0, d_{ln}(t) - r_{l} - \rho_{n})
\end{aligned}
\end{equation}
Note that $d_{ln}^{-}(t) \leq 0$ for all $t \in [0, T]$. Similar to (\ref{ZMJ_PersistentMonitoring_dpq}), a collision-free trajectory satisfies $d_{ln}^{-}(t) = 0$ for all $t$. Considering all agents and obstacles, we define the obstacle collision avoidance component $J_3(t)$ of the objective function as follows:
\begin{equation}\label{ZMJ_PersistentMonitoring_J3}
\begin{aligned}
J_{3}(t) = \displaystyle{\sum_{n = 1}^N}\displaystyle{\sum_{l = 1}^L}d_{ln}^{-}(t)
\end{aligned}
\end{equation}

Now we are ready to formulate the optimal persistent monitoring problem to obtain collision-free agent trajectories so that the cumulative uncertainty over all weighted target points $\left\{z_{1},...,z_{M}\right\}$ in the mission space is minimized over a fixed time horizon $T$.
Let $\mathbf{u}(t) = [u_{1}(t),...,u_{N}(t)]^{\top}$ and $\bm{\theta}(t) = [\theta_{1}(t),...,\theta_{N}(t)]^{\top}$. Then the optimal control problem \textbf{P1} is formulated as follows:
\begin{equation}\label{ZMJ_PersistentMonitoring_J}
\begin{aligned}
\textbf{P1: } \min \limits_{\textbf{u}(t),\bm{\theta}(t)} J = \frac{1}{T}\int_{0}^{T} \Big(J_{1}(t) + M_{2}J_{2}(t) + M_{3}J_{3}(t)\Big)dt
\end{aligned}
\end{equation}
where $M_{2}$ and $M_{3}$ are large negative numbers.

\textbf{Remark 2.} The objective function is designed on the basis of imposing collisions penalties \cite{sun2006optimization}. Along a collision-free trajectory both $J_{2}(t)$ and $J_3(t)$ are zeros. The large negative numbers $M_{2}$ and $M_{3}$ penalize all forms of collision (agent to agent, agent to obstacle) along the agent trajectory.

\section{Optimal Control Solution}\label{ZMJ_PersistentMonitoring_section_3}
In this section, we use optimal control theory \cite{bryson1975applied} to derive necessary conditions for the optimal solutions of \textbf{P1}.
We define the state vector
\begin{equation}
\begin{split}
\mathbf{y}(t) = & [s_{1}^{x}(t),s_{1}^{y}(t),...,s_{N}^{x}(t),s_{N}^{y}(t),v_{1}^{x}(t),v_{1}^{y}(t),...,\\
& v_{N}^{x}(t),v_{N}^{y}(t),R_{1}(t),...,R_{M}(t)]^{\top}
\end{split}
\end{equation}
and the associated costate vector
\begin{equation}
\begin{aligned}
\bm{\lambda}(t) = & [\lambda_{1}^{x}(t),\lambda_{1}^{y}(t),...,\lambda_{N}^{x}(t),\lambda_{N}^{y}(t),\mu_{1}^{x}(t),\mu_{1}^{y}(t),...,\\
& \mu_{N}^{x}(t),\mu_{N}^{y}(t),\gamma_{1}(t),...,\gamma_{M}(t)]^{\top}
\end{aligned}
\end{equation}
In addition, we introduce $\bm{\eta} = [\eta_{1}(t),...,\eta_{N}(t)]^{\top}$ to handle the inequality constraint on the state variable $v_{n}(t)$ such that
\begin{equation}\label{ZMJ_PersistentMonitoring_v_costate}
\begin{aligned}
\left\{
\begin{array}{rcl}
\eta_{n} = 0,\; {{\rm if}\; \Vert v_{n}(t) \Vert <  v_{n}^{max}}\\
\eta_{n} > 0,\; {{\rm if}\; \Vert v_{n}(t) \Vert =  v_{n}^{max}}
\end{array} \right.
\end{aligned}
\end{equation}
The Hamiltonian is given by
\begin{equation}\label{ZMJ_PersistentMonitoring_hamiltonian}
\begin{aligned}
H &= \Big(J_{1}(t) + M_{2}J_{2}(t) + M_{3}J_{3}(t)\Big) + \displaystyle{\sum_{i = 1}^M}\gamma_{i}(t)\dot{R}_{i}(t)\\
&+ \displaystyle{\sum_{n = 1}^N} \lambda_{n}^{x}(t)v_{n}^{x}(t) + \displaystyle{\sum_{n = 1}^N} \lambda_{n}^{y}(t)v_{n}^{y}(t)\\
&+ \displaystyle{\sum_{n = 1}^N} \mu_{n}^{x}(t)u_{n}(t)\cos\theta_{n}(t) + \displaystyle{\sum_{n = 1}^N} \mu_{n}^{y}(t)u_{n}(t)\sin\theta_{n}(t)\\
&+ \displaystyle{\sum_{n = 1}^N}\eta_{n}(t)( \Vert v_{n}(t) \Vert -  v_{n}^{max})
\end{aligned}
\end{equation}
and the costate equations can be obtained through $\dot{\bm{\lambda}} = -\frac{\partial H}{\partial \textbf{y}}$.

For the case $\mu_{n}^{y}(t) = 0$, we get
\begin{equation}\label{ZMJ_PersistentMonitoring_hamiltonian0}
\begin{aligned}
H &= \Big(J_{1}(t) + M_{2}J_{2}(t) + M_{3}J_{3}(t)\Big) + \displaystyle{\sum_{i = 1}^M}\gamma_{i}(t)\dot{R}_{i}(t)\\
&+ \displaystyle{\sum_{n = 1}^N} \lambda_{n}^{x}(t)v_{n}^{x}(t) + \displaystyle{\sum_{n = 1}^N} \lambda_{n}^{y}(t)v_{n}^{y}(t)\\
&+ \displaystyle{\sum_{n = 1}^N} \mu_{n}^{x}(t)u_{n}(t)\cos\theta_{n}(t) + \displaystyle{\sum_{n = 1}^N}\eta_{n}(t)( \Vert v_{n}(t) \Vert -  v_{n}^{max})
\end{aligned}
\end{equation}

For the case $\mu_{n}^{y}(t) \neq 0$, after some trigonometric operations on (\ref{ZMJ_PersistentMonitoring_hamiltonian}), we get
\begin{equation}
\begin{aligned}
H 
&= \Big(J_{1}(t) + M_{2}J_{2}(t) + M_{3}J_{3}(t)\Big) + \displaystyle{\sum_{i = 1}^M}\gamma_{i}(t)\dot{R}_{i}(t)\\
&+ \displaystyle{\sum_{n = 1}^N} \lambda_{n}^{x}(t)v_{n}^{x}(t) + \displaystyle{\sum_{n = 1}^N} \lambda_{n}^{y}(t)v_{n}^{y}(t)\\
&+ \displaystyle{\sum_{n = 1}^N}{\rm sgn}(\mu_{n}^{y}(t))u_{n}(t)\sqrt{(\mu_{n}^{x}(t))^{2} + (\mu_{n}^{y}(t))^{2}}\\
&\times \Bigg[ \frac{{\rm sgn}(\mu_{n}^{y}(t))\mu_{n}^{x}(t)\cos\theta_{n}(t)}{\sqrt{(\mu_{n}^{x}(t))^{2} + (\mu_{n}^{y}(t))^{2}}} + \frac{|\mu_{n}^{y}(t)|\sin\theta_{n}(t)}{\sqrt{(\mu_{n}^{x}(t))^{2} + (\mu_{n}^{y}(t))^{2}}} \Bigg]\\
&+ \displaystyle{\sum_{n = 1}^N}\eta_{n}(t)( \Vert v_{n}(t) \Vert  -  v_{n}^{max})
\end{aligned}
\end{equation}
where sgn($\cdot$) represents the sign function. Combining the trigonometric function terms, we obtain
\begin{equation}\label{ZMJ_PersistentMonitoring_hamiltonian1}
\begin{aligned}
H &= \Big(J_{1}(t) + M_{2}J_{2}(t) + M_{3}J_{3}(t)\Big) + \displaystyle{\sum_{i = 1}^M}\gamma_{i}(t)\dot{R}_{i}(t)\\
&+ \displaystyle{\sum_{n = 1}^N} \lambda_{n}^{x}(t)v_{n}^{x}(t) + \displaystyle{\sum_{n = 1}^N} \lambda_{n}^{y}(t)v_{n}^{y}(t)\\
&+ \displaystyle{\sum_{n = 1}^N}{\rm sgn}(\mu_{n}^{y}(t))u_{n}(t)\sqrt{(\mu_{n}^{x}(t))^{2} + (\mu_{n}^{y}(t))^{2}}\sin\big(\theta_{n}(t)\\
&+ \Psi_{n}(t)\big) + \displaystyle{\sum_{n = 1}^N}\eta_{n}(t)\big(\Vert v_{n}(t) \Vert -  v_{n}^{max}\big)
\end{aligned}
\end{equation}
where $\Psi_{n}(t)$ is defined so that $\tan\Psi_{n}(t) = \frac{\mu_{n}^{x}(t)}{\mu_{n}^{y}(t)}$ for $\mu_{n}^{y}(t) \neq 0$.
In the following discussion, we assume that any ``singular interval" where both $\mu_{n}^{x}(t) = 0$ and $\mu_{n}^{y}(t) = 0$ is excluded.

Note that the initial states of this problem are given but the terminal states are free. Thus, the terminal costate are $\gamma_{i}(T) = 0, i=1,...,M$, $\lambda_{n}^{x}(T) = \lambda_{n}^{y}(T) = 0$ and $\mu_{n}^{x}(T) = \mu_{n}^{y}(T) = 0, n=1,...,N$. According to \cite{bryson1975applied}, the minimal solution can be obtained by solving the Two Point Boundary Value Problem. Applying the Pontryagin Minimum Principle to (\ref{ZMJ_PersistentMonitoring_hamiltonian0}) and (\ref{ZMJ_PersistentMonitoring_hamiltonian1}) with $u_{n}^{\star}(t)$, $\theta_{n}^{\star}(t)$, $t \in [0,T)$ denoting the optimal controls, we have
\begin{equation}
\begin{aligned}
H(\textbf{y}^{\star},\bm{\lambda}^{\star},\textbf{u}^{\star},\bm{\theta}^{\star}) = \min \limits_{\textbf{u}(t),\bm{\theta}(t)} H(\textbf{y},\bm{\lambda},\textbf{u},\bm{\theta})
\end{aligned}
\end{equation}
The necessary conditions for the optimal controls are: \\
\begin{equation}
\Vert u_{n}^{\star}(t) \Vert = u_{n}^{max}
\end{equation}
and
\begin{equation}
\left\{
\begin{split}
& \sin(\theta_{n}^{\star}(t) + \Psi_{n}(t)) = 1, & \text{if}\;\; \mu_{n}^{y}(t) < 0\\
& \sin(\theta_{n}^{\star}(t) + \Psi_{n}(t)) = -1, & \text{if}\;\; \mu_{n}^{y}(t) > 0\\
& \cos\theta_{n}^{\star}(t) = 1, & \text{if}\;\; \mu_{n}^{y}(t) = 0\;\;\text{and}\;\; \mu_{n}^{x}(t) < 0\\
& \cos\theta_{n}^{\star}(t) = -1, & \text{if}\;\; \mu_{n}^{y}(t) = 0\;\;\text{and}\;\; \mu_{n}^{x}(t) > 0
\end{split}
\right.
\end{equation}

\textbf{Remark 3.} Unlike the 1D persistent monitoring analysis in \cite{wang2018optimal2} where $u_{i}^{\ast}(t)$ depends absolutely on the sign of the costate $\lambda_{u_{i}}^{\ast}(t)$, in 2D $u_{n}^{\star}(t) = u_{n}^{max}$ and the orientation $\theta_{n}^{\star}(t)$ depends on the sign of costate variables.

Note that in (\ref{ZMJ_PersistentMonitoring_v_constraint}) and (\ref{ZMJ_PersistentMonitoring_v_costate}), if $\Vert v_{n}(t) \Vert = v_{n}^{max}$, then $\dot{v}_{n}(t) = 0$. In particular, as long as the agent $n$ reaches its maximal velocity, it will keep the maximal velocity moving.

Revisiting the Hamiltonian in (\ref{ZMJ_PersistentMonitoring_hamiltonian}) and according to \cite{bryson1975applied}, $\theta_{n}^{\star}(t)$ can be obtained by solving:
\begin{equation}
\begin{aligned}
\frac{\partial H}{\partial\theta_{n}} = -\mu_{n}^{x}(t)u_{n}(t)\sin\theta_{n}(t) + \mu_{n}^{y}(t)u_{n}(t)\cos\theta_{n}(t) = 0
\end{aligned}
\end{equation}
from which we obtain
\begin{equation}
\begin{aligned}
\tan\theta_{n}^{\star}(t) = \frac{\mu_{n}^{y}(t)}{\mu_{n}^{x}(t)}
\end{aligned}
\end{equation}

So far, we know $u_{n}^{\star}(t) = u_{n}^{max}$ and there is only $\theta_{n}^{\star}(t)$ left to be evaluated. This can be accomplished by solving a standard Two Point Boundary Value Problem (TPBVP), which involves a forward integration of the states and a backward integration of the costates. However, solving the TPBVP problem is computationally intractable  as numbers of agents and targets increase. Furthermore, the switches on the state dynamics and the presence of obstacles exacerbate the computational complexity. Therefore, we will search alternatives in the next section.

\section{Agent Trajectory Parameterization And Optimization}\label{ZMJ_PersistentMonitoring_section_4}
The result of \cite{lin2015optimal} indicates that under some assumptions an elliptical trajectory outperforms a linear one when using the average uncertainty metric as a comparison criterion.
In fact, elliptical trajectories degenerate to linear ones when the minor axis of the ellipse becomes zero. Based on the result that elliptical trajectories are smooth and periodic, and are more suitable for 2D monitoring \cite{lin2015optimal}, we also select them for agents to execute the monitoring task. Under the optimal control derived in Section III, the agent first accelerates along the elliptical trajectory with the maximal acceleration $u_{n}^{max}$. Ever since it reaches the maximal velocity, it maintains the maximal velocity along the monitoring task.

\subsection{Elliptical Trajectories}\label{ZMJ_PersistentMonitoring_section_4A}
We assign each agent an elliptical trajectory parameterized by center coordinate, major axis, minor axis, and orientation and rewrite \textbf{P1} as a parametric optimization problem.

\begin{figure}
\centering
\includegraphics[height=3cm,width=6cm]{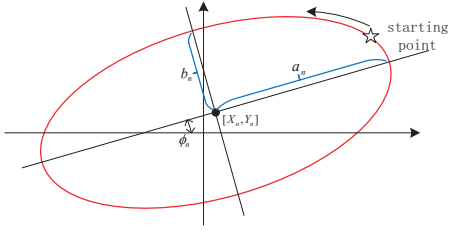}
\caption{Illustrations of the parametric form of an elliptical trajectory.} \label{ZMJ_PersistentMonitoring_tuoyuan}
\end{figure}

Define an elliptical trajectory such that the agent position $s_{n}(t)=[s_{n}^{x}(t), s_{n}^{y}(t)]$ follows the general parametric form of the ellipse:
\begin{equation}\label{ZMJ_PersistentMonitoring_s_ellipse}
\begin{aligned}
s_{n}^{x}(t) &= X_{n} + a_{n}\cos\varphi_{n}(t)\cos\phi_{n} - b_{n}\sin\varphi_{n}(t)\sin\phi_{n}\\
s_{n}^{y}(t) &= Y_{n} + a_{n}\cos\varphi_{n}(t)\sin\phi_{n} + b_{n}\sin\varphi_{n}(t)\cos\phi_{n}
\end{aligned}
\end{equation}
As illustrated in Fig. \ref{ZMJ_PersistentMonitoring_tuoyuan}, $[X_{n}, Y_{n}]$ is the center of the ellipse, $a_{n}$ and $b_{n}$ are the radii of the major and minor axes respectively, $\phi_{n} \in [0, 2\pi)$ is the orientation, and $\varphi_{n}(t)$ is the eccentric anomaly. Upon initializing its position at the starting point, the agent always moves along the ellipse.

Taking the derivatives of (\ref{ZMJ_PersistentMonitoring_s_ellipse}), we obtain the velocity information of the agent as follows,
\begin{equation}\label{ZMJ_PersistentMonitoring_v_ellipse}
\begin{aligned}
\dot{s}_{n}^{x}(t)&= -\dot{\varphi}_{n}(t)\big(a_{n}\sin\varphi_{n}(t)\cos\phi_{n} + b_{n}\cos\varphi_{n}(t)\sin\phi_{n}\big)\\
\dot{s}_{n}^{y}(t)&= -\dot{\varphi}_{n}(t)\big(a_{n}\sin\varphi_{n}(t)\sin\phi_{n} - b_{n}\cos\varphi_{n}(t)\cos\phi_{n}\big)
\end{aligned}
\end{equation}
Again, taking the derivatives of (\ref{ZMJ_PersistentMonitoring_v_ellipse}), we get the acceleration information of the agent as follows,
\begin{equation}\label{ZMJ_PersistentMonitoring_u_ellipse}
\begin{aligned}
\ddot{s}_{n}^{x}(t) =
&- a_{n}\cos\phi_{n}\big(\ddot{\varphi}_{n}(t)\sin\varphi_{n}(t)+\dot{\varphi}_{n}(t)^{2}\cos\varphi_{n}(t)\big)\\
&- b_{n}\sin\phi_{n}\big(\ddot{\varphi}_{n}(t)\cos\varphi_{n}(t)-\dot{\varphi}_{n}(t)^{2}\sin\varphi_{n}(t)\big)\\
\ddot{s}_{n}^{y}(t) =
&- a_{n}\sin\phi_{n}\big(\ddot{\varphi}_{n}(t)\sin\varphi_{n}(t)+\dot{\varphi}_{n}(t)^{2}\cos\varphi_{n}(t)\big)\\
&+ b_{n}\cos\phi_{n}\big(\ddot{\varphi}_{n}(t)\cos\varphi_{n}(t)-\dot{\varphi}_{n}(t)^{2}\sin\varphi_{n}(t)\big)
\end{aligned}
\end{equation}

Since the agent moves with constant acceleration $u_{n}^{max}$ on the elliptical trajectory (if not at the maximal velocity), we have
\begin{equation}\label{ZMJ_PersistentMonitoring_umax}
\begin{aligned}
{\ddot{s}_{n}^{x}(t)}^{2} + {\ddot{s}_{n}^{y}(t)}^{2} = \left({u_{n}^{max}}\right)^{2}
\end{aligned}
\end{equation}
\begin{equation}\label{ZMJ_PersistentMonitoring_vt}
\begin{aligned}
{\dot{s}_{n}^{x}(t)}^{2} + {\dot{s}_{n}^{y}(t)}^{2} = \Vert v_{n}(t)\Vert^{2}
\end{aligned}
\end{equation}
for some $t \in [0,T]$, where ${v_{n}(t)}$ is the current velocity of agent $n$. Applying (\ref{ZMJ_PersistentMonitoring_u_ellipse}) to (\ref{ZMJ_PersistentMonitoring_umax}) and  (\ref{ZMJ_PersistentMonitoring_v_ellipse}) to (\ref{ZMJ_PersistentMonitoring_vt}), we can obtain $\ddot{\varphi}_{n}(t)$ and $\dot{\varphi}_{n}(t)$, respectively.

Once the agent reaches the maximal velocity, it will maintain this maximal velocity thereafter and we have
\begin{equation}\label{ZMJ_PersistentMonitoring_vmax}
\begin{aligned}
{\dot{s}_{n}^{x}(t)}^{2} + {\dot{s}_{n}^{y}(t)}^{2} = \left({v_{n}^{max}}\right)^{2}
\end{aligned}
\end{equation}
for some $t \in [0,T]$. Applying (\ref{ZMJ_PersistentMonitoring_v_ellipse}) to (\ref{ZMJ_PersistentMonitoring_vmax}), we obtain $\dot{\varphi}_{n}(t)$.

Based on these processes, $\varphi_{n}(t), t \in [0,T]$, can be obtained through iterations, thus the complete trajectory of the agent on the ellipse can be obtained.
\subsection{Optimal Trajectory Design}\label{ZMJ_PersistentMonitoring_section_4B}
The trajectories of all agents can be parameterized by $\Theta = [\Theta_{1}, ..., \Theta_{N}]^{\top}$ with $\Theta_{n} = [X_{n},Y_{n},a_{n},b_{n},\phi_{n}], n = 1, ..., N$. Therefore, the objective function (\ref{ZMJ_PersistentMonitoring_J}) can be rewritten as $J(\Theta)$. We seek to obtain $\Theta^{\star} = [\Theta_{1}^{\star}, ..., \Theta_{N}^{\star}]^{\top}$ by minimizing $J(\Theta)$. We use the standard gradient descent method \cite{boyd2004convex} as follows,
\begin{equation}\label{ZMJ_PersistentMonitoring_Gradient}
\begin{aligned}
\Theta^{h} = \Theta^{h-1} - \alpha^{h}\nabla J(\Theta^{h-1})
\end{aligned}
\end{equation}
where $\alpha^{h}$ is a suitable step size and $\nabla J(\Theta^{h-1})$ is the gradient of $J$ with respect to $\Theta = [\Theta_{1}, ..., \Theta_{N}]^{\top}$. The trajectory parameters can be optimized iteratively through \eqref{ZMJ_PersistentMonitoring_Gradient} and the terminal condition is given by
\begin{equation}\label{ZMJ_PersistentMonitoring_terminal}
\begin{aligned}
|J(\Theta^{h}) - J(\Theta^{h-1})| < \varepsilon
\end{aligned}
\end{equation}
where $\varepsilon > 0$ is a prespecified constant.

\subsection{Gradient Calculation}\label{ZMJ_PersistentMonitoring_section_4C}
\textbf{\emph{1)} Infinitesimal Perturbation Analysis}. Now we briefly review Infinitesimal Perturbation Analysis (IPA) as introduced in \cite{cassandras2010perturbation}. IPA calculates the gradient of the objective function for a hybrid dynamic system which contains events leading to possible discontinuities on the gradient. Denote the continuous states of the hybrid system by $\bm{\chi}(t)$, the \emph{k}-th event time by $\rho_{k}$ and the controllable parameter by $\bm{\Theta}$. For each $[\rho_{k},\rho_{k+1})$, the state dynamics is continuous and can be written as $\dot{\bm{\chi}} = f_{k}(\bm{\chi},\bm{\Theta},t)$. The events at $\rho_{k}$ can be classified into three categories:
\begin{itemize}
\item Exogenous events. These events will cause a discrete state transition and there exists $\frac{d\rho_{k}}{d\bm{\Theta}}=0$.
\item Endogenous events. These events occur when a continuous differentiable guarding function $g_{k}(\Theta,\bm{\chi})=0$.
\item Induced events. These events are triggered by another event occurring at an earlier time.
\end{itemize}

Let $\bm{\chi} '(t) = \frac{\partial \bm{\chi}}{\partial \bm{\Theta}}$ and $\rho_{k} ' = \frac{\partial \rho_{k}}{\partial \bm{\Theta}}$. The IPA proposed in \cite{cassandras2010perturbation} shows that $\bm{\chi} '(t)$ satisfies:
\begin{equation}\label{ZMJ_PersistentMonitoring_IPA1}
\begin{aligned}
\frac{d}{dt}\bm{\chi} '(t) = \frac{\partial f_{k}}{\partial \bm{\chi}}\bm{\chi} '(t) + \frac{\partial f_{k}}{\partial \bm{\Theta}},\quad t \in [\rho_{k},\rho_{k+1})
\end{aligned}
\end{equation}
with boundary condition:
\begin{equation}\label{ZMJ_PersistentMonitoring_IPA2}
\begin{aligned}
\bm{\chi} '(\rho_{k}^{+}) = \bm{\chi} '(\rho_{k}^{-}) + [f_{k-1}(\rho_{k}^{-}) - f_{k}(\rho_{k}^{+})]\rho_{k} '
\end{aligned}
\end{equation}
If the event at $\rho_{k}$ is exogenous, then $\rho_{k} ' = 0$. Otherwise, if the event at $\rho_{k}$ is endogenous, then $\rho_{k} '$  satisfies
\begin{equation}\label{ZMJ_PersistentMonitoring_IPA3}
\begin{aligned}
\rho_{k} ' = -\Big[\frac{\partial g_{k}}{\partial \bm{\chi}}f_{k}(\rho_{k}^{-})\Big]^{-1}\Big[\frac{\partial g_{k}}{\partial \bm{\Theta}} + \frac{\partial g_{k}}{\partial \bm{\chi}}\bm{\chi} '(\rho_{k}^{-})\Big]
\end{aligned}
\end{equation}
as long as $\frac{\partial g_{k}}{\partial \bm{\chi}}f_{k}(\rho_{k}^{-}) \neq 0$ (details can be found in \cite{cassandras2010perturbation}).
Then for $t \in [\rho_{k},\rho_{k+1})$, $\bm{\chi} '(t)$ can be determined by
\begin{equation}\label{ZMJ_PersistentMonitoring_IPA4}
\begin{aligned}
\bm{\chi} '(t) = \bm{\chi} '(\rho_{k}^{+}) + \int_{\rho_{k}}^{t}\frac{d}{dt}\bm{\chi} '(t)dt
\end{aligned}
\end{equation}

Recalling Section \ref{ZMJ_PersistentMonitoring_section_2}, we summarize all events which may trigger the state transitions in our optimization problem in TABLE \ref{ZMJ_PersistentMonitoring_table_event}. The superscript $0$ indicates an event causing a variable to change from a non-zero value to zero; the superscript $+$ indicates an event causing a variable to change from zero to a positive value; and the superscript $-$ indicates an event causing a variable to change from zero to a negative value). Define $E$ as a set of ``event types" that can be associated with either an agent or a target. We use $e(\rho_{k})\in E$ to represent the event happened at time $\rho_k$.
\begin{table}[tbp]
    \setlength{\abovecaptionskip}{-0.2cm}
    \setlength{\belowcaptionskip}{-0cm}
    \renewcommand\arraystretch{1.3}
    \caption{Optimization Problem Event Set}\label{ZMJ_PersistentMonitoring_table_event}
    \begin{center}
        \begin{tabular}{|c|l|}
            \hline
            \makecell{Event Set} & \makecell{Illustration}  \\
            \hline
            $1.\; \xi_{i}^{0}$ & $R_{i}(t)$ hits $0$, for $i = 1,...,M$ \\
            \hline
            $2.\; \xi_{i}^{+}$ & $R_{i}(t)$ leaves $0$, for $i = 1,...,M$ \\
            \hline
            $3.\; u_{n}^{0}$ & $u_{n}(t)$ switches from $u_{n}^{max}$ to $0$, for $n = 1,...,N$ \\
            \hline
            $4.\; \zeta_{pq}^{0}$ & $d_{pq}^{-}(t)$ hits $0$, for $p = 1,...,N-1, q = p+1,...,N$ \\
            \hline
            $5.\; \zeta_{pq}^{-}$ & $d_{pq}^{-}(t)$ leaves $0$, for $p = 1,...,N-1, q = p+1,...,N$ \\
            \hline
            $6.\; \delta_{ln}^{0}$ & $d_{ln}^{-}(t)$ hits $0$, for $l = 1,...,L$, $n = 1,...,N$ \\
            \hline
            $7.\; \delta_{ln}^{-}$ & $d_{ln}^{-}(t)$ leaves $0$, for $l = 1,...,L$, $n = 1,...,N$ \\
            \hline
       \end{tabular}
    \end{center}
\end{table}

\textbf{\emph{2)} Gradient computation using IPA}. We write the parametric form of the objective function in (\ref{ZMJ_PersistentMonitoring_J}) as
\begin{equation}
\begin{aligned}
J(\Theta) = \frac{1}{T} \int_{0}^{T} \Big(J_{1}(\Theta) + M_{2}J_{2}(\Theta) + M_{3}J_{3}(\Theta)\Big)dt
\end{aligned}
\end{equation}
Define a time sequence $\left\{\rho_{k}(\Theta),k=1,...,K\right\}({\rm with} \rho_{0} = 0$ and $\rho_{K+1} = T)$ to describe all switching instants. Then $\nabla J(\Theta)$ can be written as
\begin{equation}\label{ZMJ_PersistentMonitoring_dJ}
\begin{aligned}
\nabla J(\Theta) &= \frac{1}{T}\displaystyle{\sum_{k = 0}^K}\Big[\int_{\rho_{k}(\Theta)}^{\rho_{k+1}(\Theta)}\Big( \nabla J_{1}(\Theta)\\
&+ M_{2}\nabla J_{2}(\Theta) + M_{3}\nabla J_{3}(\Theta)\Big)dt\Big]
\end{aligned}
\end{equation}
The evaluation of $\nabla J(\Theta)$ therefore depends entirely on $\nabla J_{1}(\Theta)$, $\nabla J_{2}(\Theta)$ and $\nabla J_{3}(\Theta)$. Note that there may be discontinuities in these derivatives and the effects of such discontinuities can be captured by IPA. Next we will focus on seeking $\nabla J_{1}(\Theta_{n})$, $\nabla J_{2}(\Theta_{n})$ and $\nabla J_{3}(\Theta_{n})$ for every agent $n$. Referring to
(\ref{ZMJ_PersistentMonitoring_J1}), (\ref{ZMJ_PersistentMonitoring_J2}) and (\ref{ZMJ_PersistentMonitoring_J3}), we can get
\begin{equation}\label{ZMJ_PersistentMonitoring_dJ1}
\begin{aligned}
\nabla J_{1}(\Theta_{n}) = \displaystyle{\sum_{i = 1}^M}\sigma_{i}\nabla R_{i}(t)
\end{aligned}
\end{equation}
\begin{equation}\label{ZMJ_PersistentMonitoring_dJ2}
\begin{aligned}
\nabla J_{2}(\Theta_{n}) =\displaystyle{\sum_{p \neq n, p \in \{1,\ldots, N\} }} \nabla d_{pn}^{-}(t)
\end{aligned}
\end{equation}
\begin{equation}\label{ZMJ_PersistentMonitoring_dJ3}
\begin{aligned}
\nabla J_{3 }(\Theta_{n}) = \displaystyle{\sum_{l = 1}^L} \nabla d_{ln}^{-}(t)
\end{aligned}
\end{equation}
We need to calculate $\nabla R_{i}(t) = [\frac{\partial R_{i}(t)}{\partial X_{n}},\frac{\partial R_{i}(t)}{\partial Y_{n}}, \frac{\partial R_{i}(t)}{\partial a_{n}}, \frac{\partial R_{i}(t)}{\partial b_{n}}, \\\frac{\partial R_{i}(t)}{\partial \phi_{n}}]^{\top}$, $\nabla d_{pn}^{-}(t) = [\frac{\partial d_{pn}^{-}(t)}{\partial X_{n}},\frac{\partial d_{pn}^{-}(t)}{\partial Y_{n}}, \frac{\partial d_{pn}^{-}(t)}{\partial a_{n}}, \frac{\partial d_{pn}^{-}(t)}{\partial b_{n}}, \\\frac{\partial d_{pn}^{-}(t)}{\partial \phi_{n}}]^{\top}$ and $\nabla d_{ln}^{-}(t) = [\frac{\partial d_{ln}^{-}(t)}{\partial X_{n}},\frac{\partial d_{ln}^{-}(t)}{\partial Y_{n}}, \frac{\partial d_{ln}^{-}(t)}{\partial a_{n}}, \frac{\partial d_{ln}^{-}(t)}{\partial b_{n}}, \\\frac{\partial d_{ln}^{-}(t)}{\partial \phi_{n}}]^{\top}$.

\textbf{\emph{2.1)} Calculate $\nabla R_{i}(t)$}. From (\ref{ZMJ_PersistentMonitoring_Uncertainty}) and (\ref{ZMJ_PersistentMonitoring_IPA4}), for $t \in [\rho_{k},\rho_{k+1})$, we can obtain
\begin{equation}\label{ZMJ_PersistentMonitoring_R_IPA4}
\begin{aligned}
\nabla R_{i}(t) = \nabla R_{i}(\rho_{k}^{+})
+ \left\{
\begin{array}{rcl}
0,\;{{\rm if}\; R_{i}(t) = 0, A_{i} < BP_{i}(\mathbf{x}(t))}\\
\int_{\rho_{k}}^{t}\frac{d}{dt}\nabla R_{i}(t)dt,\qquad\;{{\rm otherwise}}
\end{array} \right.
\end{aligned}
\end{equation}
where the integral term calculated by (\ref{ZMJ_PersistentMonitoring_IPA1}) gives
\begin{equation}
\begin{aligned}
\frac{d}{dt}\nabla R_{i}(t) = - B\nabla P_{i}(\mathbf{x}(t))
\end{aligned}
\end{equation}
From (\ref{ZMJ_PersistentMonitoring_probability}), we can get
\begin{equation}\label{ZMJ_PersistentMonitoring_dpro}
\begin{aligned}
\nabla P_{i}(\mathbf{x}(t)) &= \frac{\partial P_{i}(\mathbf{x}(t))}{\partial s_{n}^{x}(t)}\nabla s_{n}^{x}(t) + \frac{\partial P_{i}(\mathbf{x}(t))}{\partial s_{n}^{y}(t)}\nabla s_{n}^{y}(t)\\
&+ \frac{\partial P_{i}(\mathbf{x}(t))}{\partial v_{n}^{x}(t)}\nabla v_{n}^{x}(t) + \frac{\partial P_{i}(\mathbf{x}(t))}{\partial v_{n}^{y}(t)}\nabla v_{n}^{y}(t)
\end{aligned}
\end{equation}
where the derivative terms in (\ref{ZMJ_PersistentMonitoring_dpro}) can be easily obtained from (\ref{ZMJ_PersistentMonitoring_probability}). Then, $\nabla s_{n}^{x}(t)$, $\nabla s_{n}^{y}(t)$, $\nabla v_{n}^{x}(t)$ and $\nabla v_{n}^{y}(t)$ can be obtained from (\ref{ZMJ_PersistentMonitoring_s_ellipse}) and (\ref{ZMJ_PersistentMonitoring_v_ellipse}) as follows,
\begin{equation}\label{ZMJ_PersistentMonitoring_sx_ellipse}
\begin{aligned}
\frac{\partial s_{n}^{x}(t)}{\partial X_{n}} &= 1,\quad \frac{\partial s_{n}^{x}(t)}{\partial a_{n}} = \cos\varphi_{n}(t)\cos\phi_{n}\\
\frac{\partial s_{n}^{x}(t)}{\partial Y_{n}} &= 0,\quad \frac{\partial s_{n}^{x}(t)}{\partial b_{n}} = -\sin\varphi_{n}(t)\sin\phi_{n}\\
\frac{\partial s_{n}^{x}(t)}{\partial \phi_{n}} &= - a_{n}\cos\varphi_{n}(t)\sin\phi_{n} - b_{n}\sin\varphi_{n}(t)\cos\phi_{n}
\end{aligned}
\end{equation}
\begin{equation}\label{ZMJ_PersistentMonitoring_sy_ellipse}
\begin{aligned}
\frac{\partial s_{n}^{y}(t)}{\partial X_{n}} &= 0,\quad \frac{\partial s_{n}^{y}(t)}{\partial a_{n}} = \cos\varphi_{n}(t)\sin\phi_{n}\\
\frac{\partial s_{n}^{y}(t)}{\partial Y_{n}} &= 1,\quad \frac{\partial s_{n}^{y}(t)}{\partial b_{n}} = \sin\varphi_{n}(t)\cos\phi_{n}\\
\frac{\partial s_{n}^{y}(t)}{\partial \phi_{n}} &= a_{n}\cos\varphi_{n}(t)\cos\phi_{n} - b_{n}\sin\varphi_{n}(t)\sin\phi_{n}
\end{aligned}
\end{equation}
\begin{equation}\label{ZMJ_PersistentMonitoring_vx_ellipse}
\begin{aligned}
\frac{\partial v_{n}^{x}(t)}{\partial X_{n}} &= 0,\quad \frac{\partial v_{n}^{x}(t)}{\partial a_{n}} = -\dot{\varphi}_{n}(t)\sin\varphi_{n}(t)\cos\phi_{n}\\
\frac{\partial v_{n}^{x}(t)}{\partial Y_{n}} &= 0,\quad \frac{\partial v_{n}^{x}(t)}{\partial b_{n}} = -\dot{\varphi}_{n}(t)\cos\varphi_{n}(t)\sin\phi_{n}\\
\frac{\partial v_{n}^{x}(t)}{\partial \phi_{n}} &= \dot{\varphi}_{n}(t)\big(a_{n}\sin\varphi_{n}(t)\sin\phi_{n} - b_{n}\cos\varphi_{n}(t)\cos\phi_{n}\big)
\end{aligned}
\end{equation}
\begin{equation}\label{ZMJ_PersistentMonitoring_vy_ellipse}
\begin{aligned}
\frac{\partial v_{n}^{y}(t)}{\partial X_{n}} &= 0,\quad \frac{\partial v_{n}^{y}(t)}{\partial a_{n}} = - \dot{\varphi}_{n}(t)\sin\varphi_{n}(t)\sin\phi_{n}\\
\frac{\partial v_{n}^{y}(t)}{\partial Y_{n}} &= 0,\quad \frac{\partial v_{n}^{y}(t)}{\partial b_{n}} = \dot{\varphi}_{n}(t)\cos\varphi_{n}(t)\cos\phi_{n}\\
\frac{\partial v_{n}^{y}(t)}{\partial \phi_{n}} &= -\dot{\varphi}_{n}(t)\big(a_{n}\sin\varphi_{n}(t)\cos\phi_{n} + b_{n}\cos\varphi_{n}(t)\sin\phi_{n}\big)
\end{aligned}
\end{equation}

\textbf{Remark 4.} Note that the calculations of (\ref{ZMJ_PersistentMonitoring_sx_ellipse})-(\ref{ZMJ_PersistentMonitoring_vy_ellipse}) depend on $\varphi_{n}(t)$ and $\dot{\varphi}_{n}(t)$.
When the agent $n$ is not at the maximal velocity, $\varphi_{n}(t)$ and $\dot{\varphi}_{n}(t)$ can be obtained by solving (\ref{ZMJ_PersistentMonitoring_umax}) and (\ref{ZMJ_PersistentMonitoring_vt}). When the agent $n$ is at the maximal velocity, i.e. the event $\left\{u_{n}^{0} \right\}$ in TABLE \ref{ZMJ_PersistentMonitoring_table_event} occurs, they can be obtained by solving (\ref{ZMJ_PersistentMonitoring_vmax}).

The following proposition derives the value of $\nabla R_{i}(\rho_{k}^{+})$ in (\ref{ZMJ_PersistentMonitoring_R_IPA4}) after the event time $t = \rho_{k}$. Note that the calculation of $\nabla R_{i}(\rho_{k}^{+})$ only involves the first two events in TABLE \ref{ZMJ_PersistentMonitoring_table_event}.

\textbf{Proposition 1.} If an event occurs at $t = \rho_{k}$, the state derivative $\nabla R_{i}(\rho_{k}^{+})$ satisfies:\\
\begin{equation}\label{ZMJ_PersistentMonitoring_GradientUncertainty}
\begin{aligned}
\nabla R_{i}(\rho_{k}^{+}) =  \left\{
\begin{array}{rcl}
0,\quad&{{{\rm if}\; e(\rho_{k}) = \xi_{i}^{0}}}\\
\nabla R_{i}(\rho_{k}^{-}),\quad&{{\rm if}\; e(\rho_{k}) = \xi_{i}^{+}}
\end{array} \right.
\end{aligned}
\end{equation}

\textbf{Proof.} See Appendix \ref{ZMJ_PersistentMonitoring_appendixA}.
$\hfill\blacksquare$

\textbf{\emph{2.2)} Calculate $\nabla d_{pn}^{-}(t)$}. This involves the derivative of the min function in (\ref{ZMJ_PersistentMonitoring_dpq}). For a given pair of agents $p$ and $n$, following the definition of $d_{pn}^{-}(t)$ in \eqref{ZMJ_PersistentMonitoring_dpq},
\begin{equation}
\begin{aligned}
\nabla d_{pn}^{-}(t) =  \left\{
\begin{split}
& 0, & \text{ if }{t \in [\zeta_{pn}^{0}, \zeta_{pn}^{-}]}\;\\
&\nabla \| s_{p}(t) - s_{n}(t) \|, & \text{ if }{t \in (\zeta_{pn}^{-}, \zeta_{pn}^{0})}
\end{split} \right.
\end{aligned}
\end{equation}

\textbf{\emph{2.3)} Calculate $\nabla d_{ln}^{-}(t)$}. This involves the derivative of the min function (\ref{ZMJ_PersistentMonitoring_dln}). For a given pair $(l,n)$ defined by obstacle $l$ and agent $n$, $\nabla d_{ln}^{-}(t)$ can be calculated as follows:
\begin{equation}
\begin{aligned}
\nabla d_{ln}^{-}(t) =  \left\{
\begin{split}
& 0, & \text{ if } t \in [\delta_{ln}^{0}, \delta_{ln}^{-}]\;\\
& \nabla \| \omega_{l}(t) - s_{n}(t) \|, & \text{ if } t \in (\delta_{ln}^{-}, \delta_{ln}^{0})
\end{split} \right.
\end{aligned}
\end{equation}

\begin{figure}[!t]
\centering
\includegraphics[height=5.4cm,width=6.7cm]{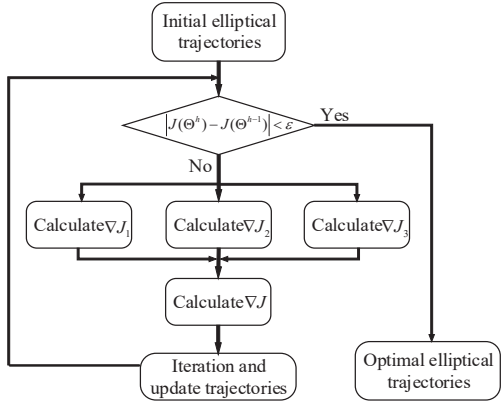}
\caption{Flow chart of Algorithm 1.}
\label{ZMJ_PersistentMonitoring_Flowchart}
\end{figure}

Based on the above analysis, the procedure for obtaining optimal elliptical trajectories is summarized in Algorithm \ref{ZMJ_PersistentMonitoring_algorithm1}, which is an off-line centralized algorithm. From Algorithm \ref{ZMJ_PersistentMonitoring_algorithm1}, we can get optimal  parameters for every agent. Then each agent will execute the monitoring task in an open-loop manner, which has been mentioned in \textbf{Remark 1}. In addition, the flow chart of the algorithm is presented in Fig. \ref{ZMJ_PersistentMonitoring_Flowchart}.

\begin{algorithm}[!t]\label{ZMJ_PersistentMonitoring_algorithm1}
    \caption{IPA-based iteration algorithm}
    \KwIn{The initial parametric trajectory $\Theta_{n}^{1}=[X_{n}^{1}, Y_{n}^{1}, a_{n}^{1}, b_{n}^{1}, \phi_{n}^{1}]$, maximal acceleration $u_{n}^{max}$, maximal velocity $v_{n}^{max}$, sensing range $r_{n}$, velocity threshold $\beta_{n}$ for effective monitoring, safety radius $\rho_{n}$, initial values of uncertainty for all targets $R_{i}(0)$, positions of obstacles,  terminal condition $\varepsilon$, iterative index $h = 1$, and $J(\Theta^{0}) = 0$.}
    \KwOut{Optimal parameters $[X_{n}^{\star}, Y_{n}^{\star}, a_{n}^{\star}, b_{n}^{\star}, \phi_{n}^{\star}]$.}
    \While{$|J(\Theta^{h}) - J(\Theta^{h-1})| > \varepsilon$}{
        Compute agent trajectory $s_{n}(t)$, $t \in [0, T]$ defined by $[X_{n}^{h}, Y_{n}^{h}, a_{n}^{h}, b_{n}^{h}, \phi_{n}^{h}]$ through (\ref{ZMJ_PersistentMonitoring_s_ellipse})-(\ref{ZMJ_PersistentMonitoring_vmax})\;
        \For{$i=1:M$}
        {
        Calculate $\dot R_{i}(t)$ according to (\ref{ZMJ_PersistentMonitoring_Uncertainty})\;
        \eIf{$\dot R_{i}(t)$ switches}
        {
            Compute $\nabla R_{i}(t)$ according to (\ref{ZMJ_PersistentMonitoring_GradientUncertainty})\;
        }
        {
            Compute $\nabla R_{i}(t)$ according to (\ref{ZMJ_PersistentMonitoring_R_IPA4})\;
        }
        }
        Compute $\nabla J_{1}$ according to (\ref{ZMJ_PersistentMonitoring_dJ1})\;
        \For{agent $p \neq$ $n$}
        {
        Calculate $d_{pn}^{-}(t)$ according to (\ref{ZMJ_PersistentMonitoring_dpq})\;
        \eIf{$d_{pn}^{-}(t) = d_{pn}(t) - \rho_{p} - \rho_{n}$}
        {
            Calculate $\nabla d_{pn}^{-}(t)$ using $\nabla s_{n}(t)$ through (\ref{ZMJ_PersistentMonitoring_sx_ellipse}) and (\ref{ZMJ_PersistentMonitoring_sy_ellipse})\;
        }
        {
            $\nabla d_{pn}^{-}(t) = 0$\;
        }
        }
        Compute $\nabla J_{2}$ according to (\ref{ZMJ_PersistentMonitoring_dJ2})\;
        \For{$l=1:L$}
        {
        Calculate $d_{ln}^{-}(t)$ according to (\ref{ZMJ_PersistentMonitoring_dln})\;
        \eIf{$d_{ln}^{-}(t) = d_{ln}(t) - r_{l} - \rho_{n}$}
        {
            Calculate $\nabla d_{ln}^{-}(t)$ using $\nabla s_{n}(t)$ through (\ref{ZMJ_PersistentMonitoring_sx_ellipse}) and (\ref{ZMJ_PersistentMonitoring_sy_ellipse})\;
        }
        {
            $\nabla d_{ln}^{-}(t) = 0$\;
        }
        }
        Compute $\nabla J_{3}$ according to (\ref{ZMJ_PersistentMonitoring_dJ3})\;
        Compute the overall gradient $\nabla J$ using (\ref{ZMJ_PersistentMonitoring_dJ})\;
        Update $[X_{n}^{h}, Y_{n}^{h}, a_{n}^{h}, b_{n}^{h}, \phi_{n}^{h}]$ through (\ref{ZMJ_PersistentMonitoring_Gradient})\; Set $h = h+1$\;
        }
        return $[X_{n}^{\star}, Y_{n}^{\star}, a_{n}^{\star}, b_{n}^{\star}, \phi_{n}^{\star}]=[X_{n}^{h}, Y_{n}^{h}, a_{n}^{h}, b_{n}^{h}, \phi_{n}^{h}]$\;
\end{algorithm}

\section{Simulation Examples}\label{ZMJ_PersistentMonitoring_section_5}
In this section, simulations are provided to illustrate the results obtained by Algorithm \ref{ZMJ_PersistentMonitoring_algorithm1}.

In the following examples, the mission space is $\bm{\Omega} \equiv [0, 10] \times [0, 5] \subset \mathbb{R}^{2}$. Note that in our work, the event excitation problem discussed in \cite{zhou2018optimal} may be caused by two situations including initial agent trajectories not visiting any targets and targets with large weights being not monitored. In either case, the performance metric gradients may not be excited, resulting in zero values that prevent the iterative algorithm (\ref{ZMJ_PersistentMonitoring_Gradient}) from converging to an optimum. This problem can be overcome by introducing a ``potential field" as detailed in earlier work \cite{zhou2018optimal}. Here, to avoid the occurrence of event excitation problem, select $\left\{z_{i} = [z_{i}^{x}, z_{i}^{y}], z_{i}^{x} = 0,1,..,10, z_{i}^{y} = 0,1,..,5\right\}$ as the target points, which are uniformly distributed in the mission space. The initial values of uncertainty for these points are $R_{i}(0) = 0, i = 1, ..., 66$. In addition, the effective sensing range for all agents is $r_{n} = 2$, the velocity threshold for effective monitoring is $\beta_{n} = 5$, the maximal acceleration is $u_{n}^{max} = 1$ and the maximal velocity is $v_{n}^{max} = 1.5$. Moreover, the safety radius $\rho_{n}$ is 0.2 for any agent $n$. Further, for avoiding slight collisions, i.e. $M_{2}J_{2}(t)$ and $M_{3}J_{3}(t)$ ending up as small numbers, we set an extra safety distance as 0.02 in the simulation. $M_{2}$ and $M_{3}$ are selected as -30000.


\textbf{\emph{Case A.} One agent case.}
In this one agent case, we use both the elliptical trajectory and the Fourier trajectory (Detailed derivation results are shown in Appendix \ref{ZMJ_PersistentMonitoring_appendixB}) to illustrate the effectiveness of Algorithm \ref{ZMJ_PersistentMonitoring_algorithm1} respectively. Since the objective function $J$ is non-convex, there may be many local optima depending on initial trajectories. Therefore, we apply multiple initial parameters and select the best. Then on the basis of elliptical trajectory, we further give some examples of comparison with two cases, i.e. monitoring with no obstacles and monitoring with the existing sensing model in \cite{lin2015optimal}, to illustrate the effectiveness of our collision avoidance method and of the new sensing model. The time horizon is set as $40s$. Let all target points have equal importance such that $\sigma_{i} = 1, i = 1, ..., 66$. The increasing rate is $A_{i} = 1, i = 1, ..., 66$ and the decreasing rate is $B = 15$.

From (\ref{ZMJ_PersistentMonitoring_s_ellipse}), the positions of agent $n$ come down to the calculation of $\varphi_{n}(t)$. Further, since the ellipse is a curve, we take the discrete way to calculate the positions. Denote $\Delta t$ as a small time interval. When the agent $n$ is not at the maximal velocity, $\varphi_{n}(t)$ can be determined by repeatedly calculating $\varphi_{n}(m+1) = \varphi_{n}(m) + \dot{\varphi}_{n}(m)\Delta t + \frac{1}{2}\ddot{\varphi}_{n}(m)\Delta t^{2}$, where $\dot{\varphi}_{n}(m)$ can be obtained by solving ${\dot{s}_{n}^{x}(t)}^{2} + {\dot{s}_{n}^{y}(t)}^{2} = {v_{n}(t)}^{2}$ ($v_{n}(t)$ is the current velocity) and $\ddot{\varphi}_{n}(m)$ can be obtained by solving (\ref{ZMJ_PersistentMonitoring_umax}). When the agent $n$ is at the maximal velocity, i.e. the event $\left\{u_{n}^{0} \right\}$ in TABLE \ref{ZMJ_PersistentMonitoring_table_event} occurs, $\varphi_{n}(m+1) = \varphi_{n}(m) + \dot{\varphi}_{n}(m)\Delta t$, where $\dot{\varphi}_{n}(m)$ can be obtained by solving (\ref{ZMJ_PersistentMonitoring_vmax}). Through all these processes, we can obtain $\varphi_{n}(t)$, $\dot{\varphi}_{n}(t)$ and $s_{n}(t)$ for all $t$. Therefore, the derivatives from (\ref{ZMJ_PersistentMonitoring_sx_ellipse}) to (\ref{ZMJ_PersistentMonitoring_vy_ellipse}) can be calculated.

\textbf{\emph{Example 1.} Illustration example of Algorithm \ref{ZMJ_PersistentMonitoring_algorithm1}.}

\begin{figure}[!h]
\centering
\subfigure[Green ellipse: initial trajectory. Blue ellipse: final trajectory. Yellow circular areas: areas covering obstacles. Red $*$: target points. Black pentagrams: starting points.]{\label{ZMJ_PersistentMonitoring_fig.41} 
\includegraphics[height=3.5cm,width=7cm]{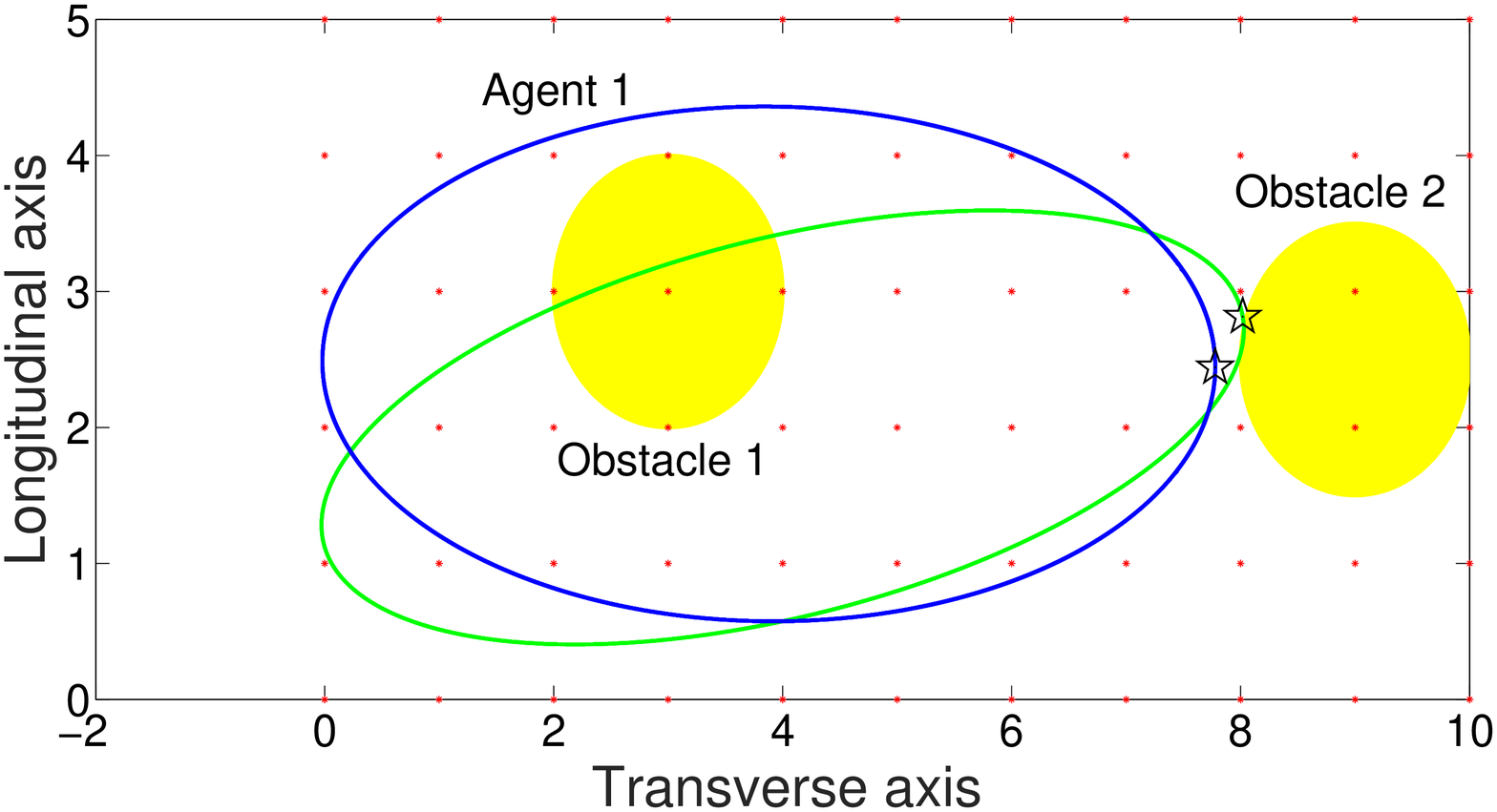}}
\subfigure[The evolution of performance metric J.]{\label{ZMJ_PersistentMonitoring_fig.42} 
\includegraphics[height=3.5cm,width=7cm]{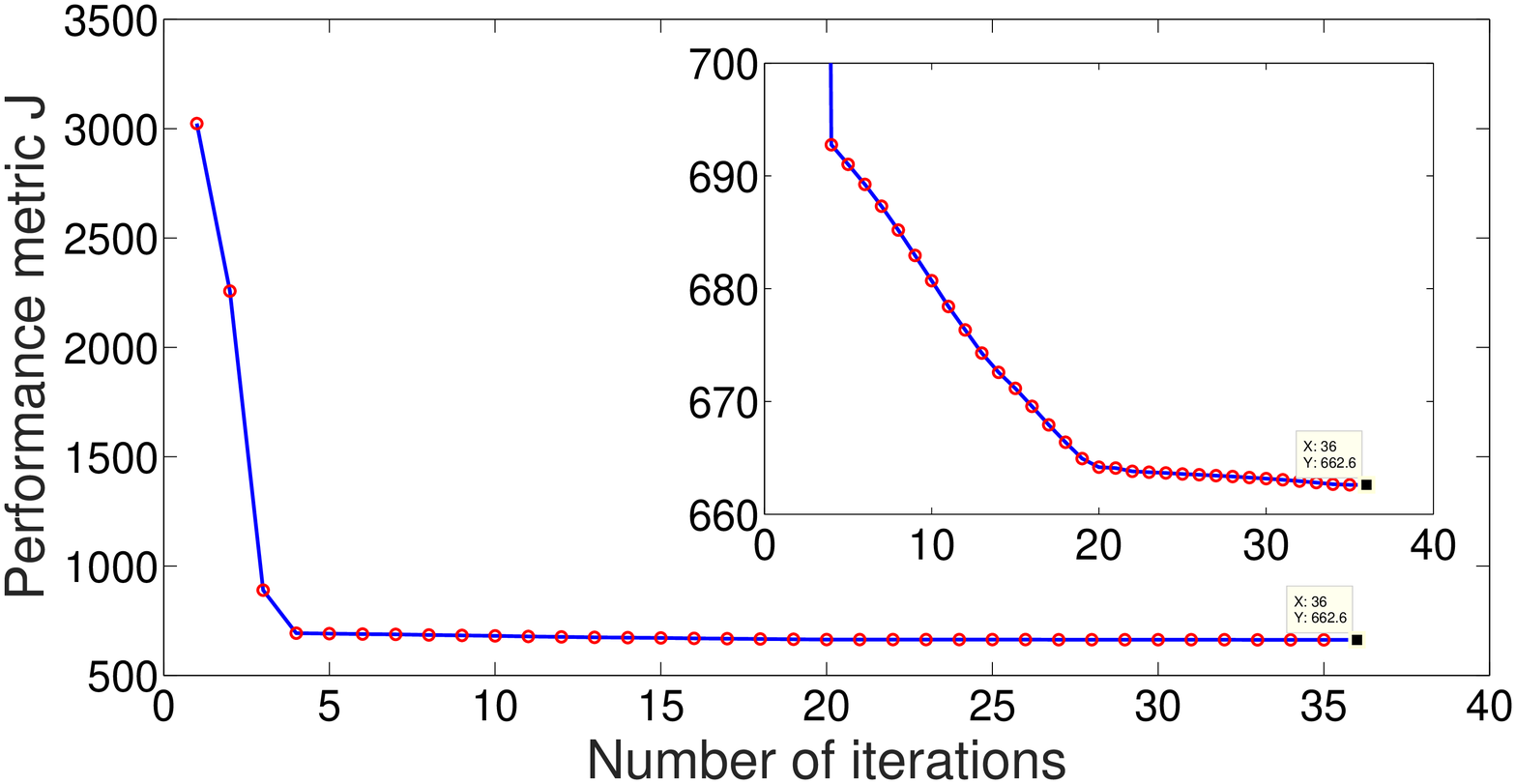}}
\subfigure[Distances between Agent 1 and Obstacles.]{\label{ZMJ_PersistentMonitoring_fig.44} 
\includegraphics[height=3.5cm,width=7cm]{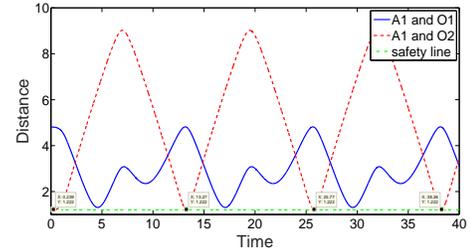}}
\subfigure[Velocity component of this work.]{\label{ZMJ_PersistentMonitoring_fig.45} 
\includegraphics[height=3.5cm,width=7cm]{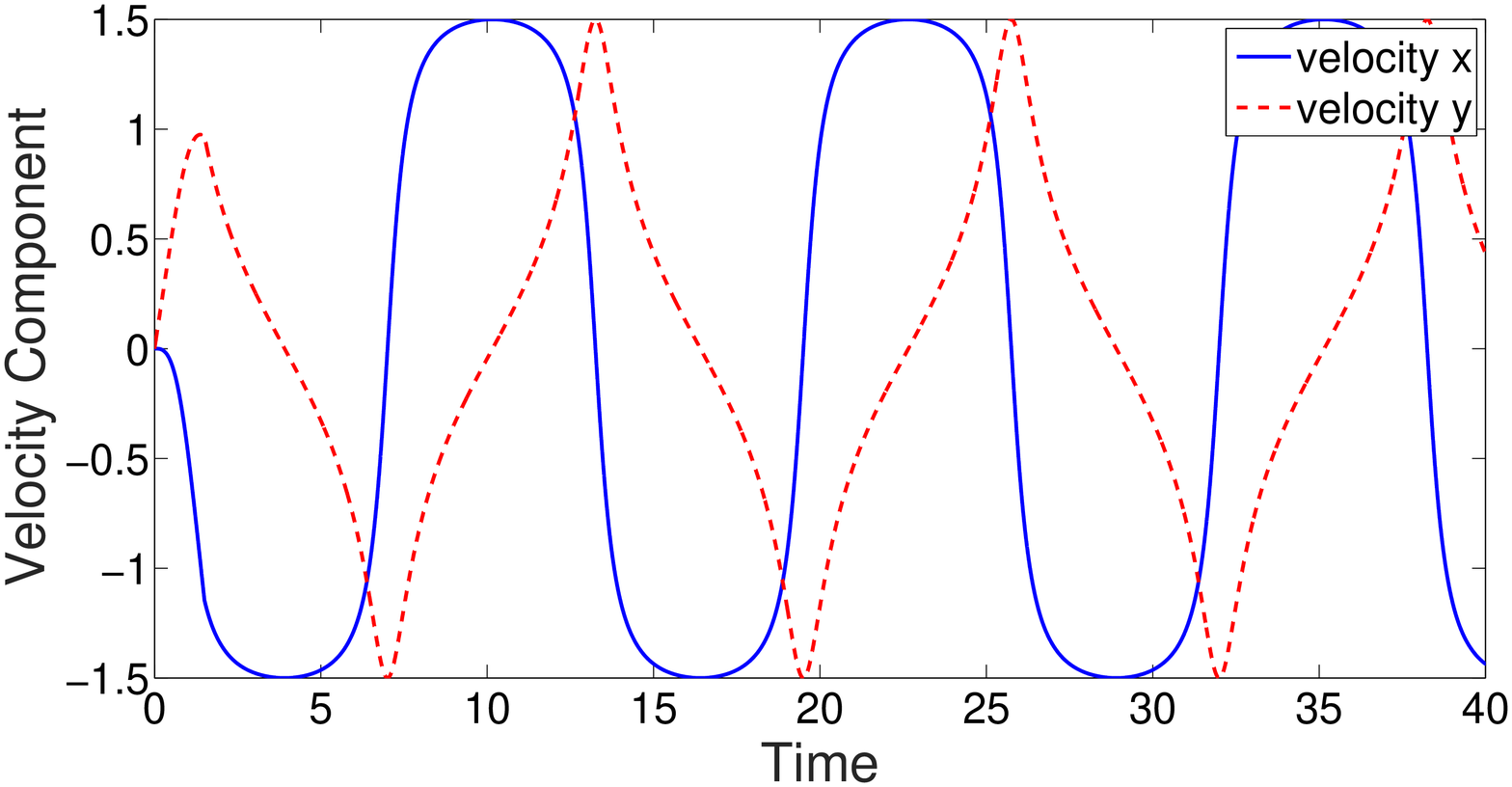}}
\subfigure[Velocity component of \cite{lin2015optimal}.]{\label{ZMJ_PersistentMonitoring_fig.46} 
\includegraphics[height=3.5cm,width=7cm]{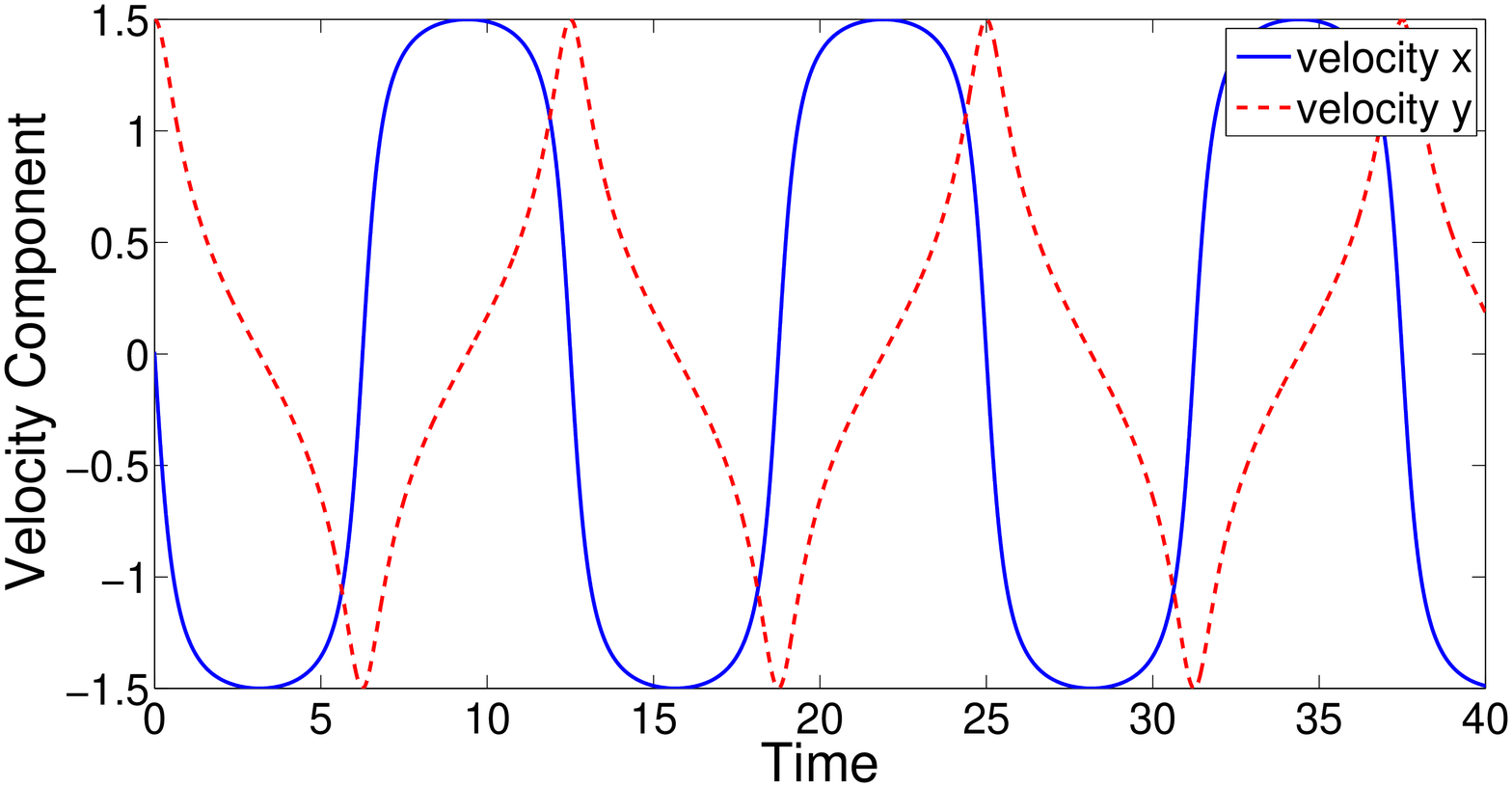}}
\caption{Elliptical trajectory: persistent monitoring task with obstacles using the IPA-based iteration algorithm for one agent.} \label{ZMJ_PersistentMonitoring_fig.4}
\end{figure}

In Fig. \ref{ZMJ_PersistentMonitoring_fig.4}, there is a persistent monitoring task with obstacles executed by one agent moving on an elliptical trajectory and we show the best from multiple elliptical results. In Fig. \ref{ZMJ_PersistentMonitoring_fig.41}, the mission space is given where there exist obstacles covered by yellow circular areas whose centers are $[3, 3]$, $[9, 2.5]$, respectively and radii are both $1$. Agent 1 moves counterclockwise on the trajectory with the pentagram as the starting point (subsequent results are the same).
Fig. \ref{ZMJ_PersistentMonitoring_fig.42} indicates that the performance metric decreases as the number of iterations increases and ultimately converges, which verifies the effectiveness of our IPA-based iteration algorithm. As we can see, there is a very abrupt drop in the first several iterations and after it the performance metric decreases with a relatively slow rate, which indicates that the initial trajectory with obstacle collisions has a very large cost and Agent 1 will continue to search for a better monitoring trajectory after avoiding obstacles.
The final performance metric $J(\Theta^{36}) = 662.6$ and $|J(\Theta^{36}) - J(\Theta^{35})| < \varepsilon = 0.01$.
According to Fig. \ref{ZMJ_PersistentMonitoring_fig.44}, we can intuitively observe that distances between Agent 1 and obstacles are greater than 1.2 all the time, which indicates that there are no collisions in the monitoring task. Actually, all simulaion results ultimately have no collisions, which shows the correctness of our objective function setting.
Furthermore, detailed velocity components of this work and of \cite{lin2015optimal} are shown in Fig. \ref{ZMJ_PersistentMonitoring_fig.45} and Fig. \ref{ZMJ_PersistentMonitoring_fig.46}, respectively. Compared with Fig. \ref{ZMJ_PersistentMonitoring_fig.46} where the agent starts with a certain velocity, the result in Fig. \ref{ZMJ_PersistentMonitoring_fig.45} can better mimic agent behavior in practice.

With the same setting as the elliptical example, simulation results of the Fourier trajectory example with one agent can be seen in Appendix \ref{ZMJ_PersistentMonitoring_appendixC}. Furthermore, the final numerical results of the two trajectories are presented in TABLE \ref{ZMJ_PersistentMonitoring_table2}, and the final costs of collisions (agent to agent, agent to obstacle), i.e. $J_{2}^{\star}$ and $J_{3}^{\star}$, are zero, which indicates that we achieve collision-free monitoring.
\begin{table}[!h]
    \setlength{\abovecaptionskip}{-0.2cm}
    \setlength{\belowcaptionskip}{-0cm}
    \renewcommand\arraystretch{1.3}
    \caption{Result Comparison For Case A}\label{ZMJ_PersistentMonitoring_table2}
    \begin{center}
        \begin{tabular}{|c|c|c|c|}
            \hline
            $Method$ & $J^{\star}$ & $J_{2}^{\star}$ & $J_{3}^{\star}$ \\
            \hline
            $Elliptical$ & $662.6$ & $0$ & $0$ \\
            \hline
            $Fourier$ & $654$ & $0$ & $0$ \\
            \hline
       \end{tabular}
    \end{center}
\end{table}
\begin{figure}[!t]
\centering
\subfigure[Green ellipse: initial trajectory. Blue ellipse: final trajectory. Red $*$: target points. Black pentagrams: starting points.]{\label{ZMJ_PersistentMonitoring_fig.51} 
\includegraphics[height=3.5cm,width=7cm]{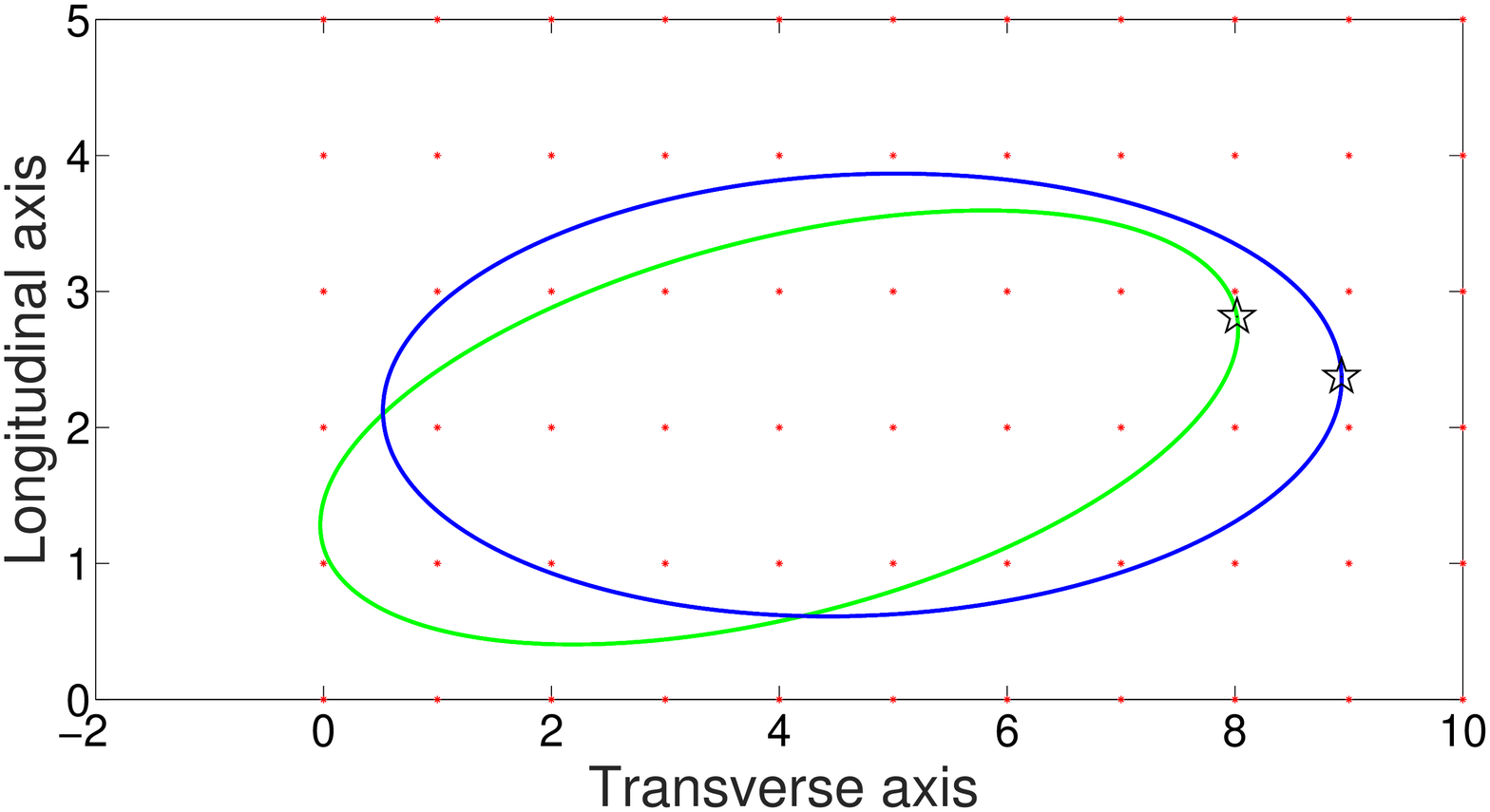}}
\subfigure[The evolution of performance metric J.]{\label{ZMJ_PersistentMonitoring_fig.52} 
\includegraphics[height=3.5cm,width=7cm]{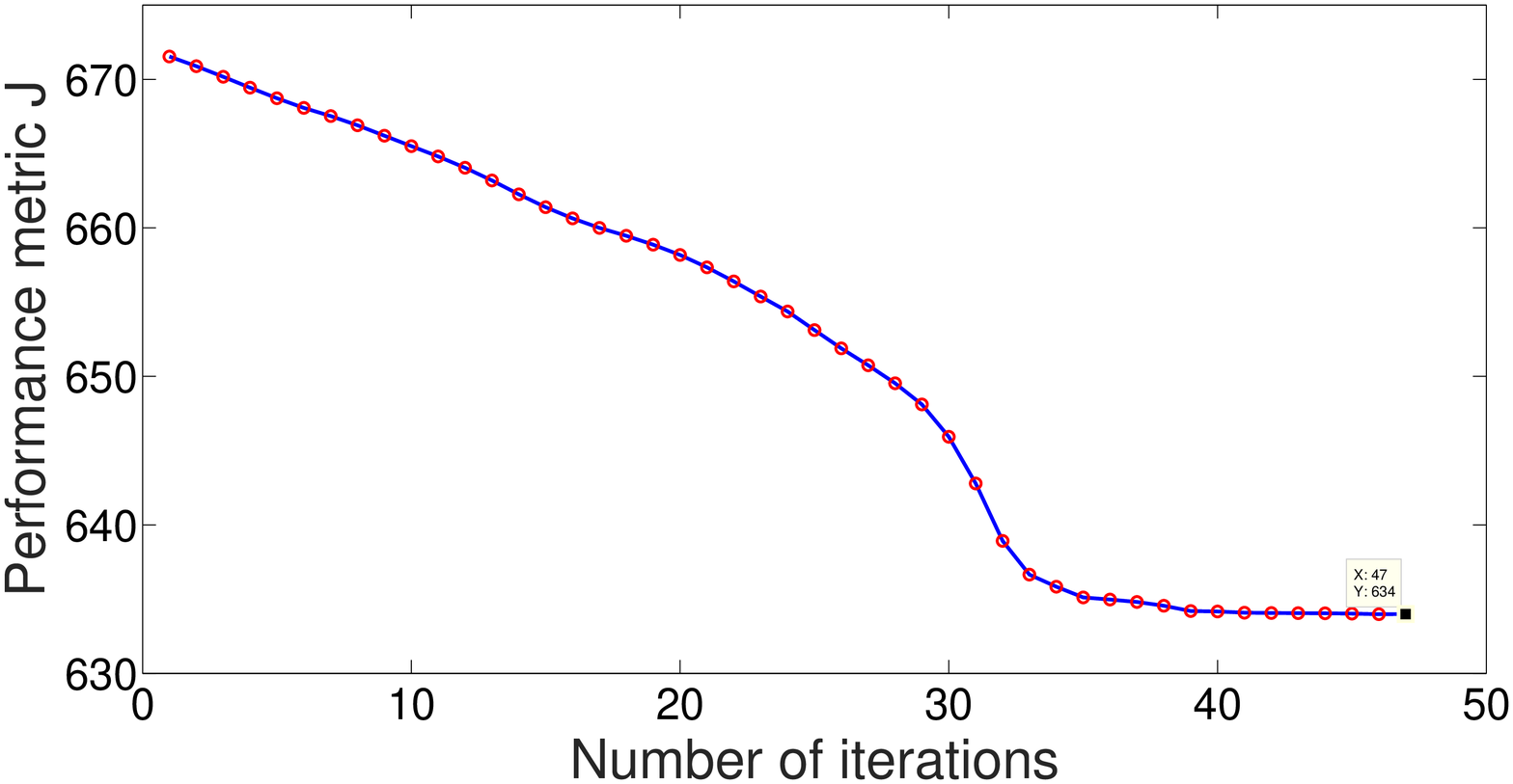}}
\caption{Elliptical trajectory: persistent monitoring task without obstacles using the IPA-based iteration algorithm for one agent.} \label{ZMJ_PersistentMonitoring_fig.5}
\end{figure}

\textbf{\emph{Example 2.} Comparison with the case without obstacles.}

In this example, we set an idealized mission space like \cite{lin2015optimal} for agents, namely, without obstacles. With the same initial conditions as in Fig. \ref{ZMJ_PersistentMonitoring_fig.4} of \textbf{\emph{Example 1}} except that there are no obstacles, the results are shown in Fig. \ref{ZMJ_PersistentMonitoring_fig.5}. Fig. \ref{ZMJ_PersistentMonitoring_fig.51} presents the optimal elliptical trajectory obtained by Algorithm \ref{ZMJ_PersistentMonitoring_algorithm1}. Due to the potential slight collisions, it is unreasonable to use the final trajectory in Fig. \ref{ZMJ_PersistentMonitoring_fig.51} to execute the persistent monitoring task with obstacles. From Fig. \ref{ZMJ_PersistentMonitoring_fig.52}, the performance metric converges to $J(\Theta^{47}) = 634$ and $|J(\Theta^{47}) - J(\Theta^{46})| < \varepsilon = 0.01$. Compared with the performance metric in Fig. \ref{ZMJ_PersistentMonitoring_fig.42}, our method can avoid collisions at the expense of some monitoring, which is quite practical in monitoring tasks where collisions with obstacles are strictly prohibited.
\begin{figure}[!h]
\centering
\subfigure[Green ellipse: initial trajectory. Blue ellipse: final trajectory. Yellow circular areas: areas covering obstacles. Red $*$: target points. Black pentagrams: starting points.]{\label{ZMJ_PersistentMonitoring_fig.61} 
\includegraphics[height=3.5cm,width=7cm]{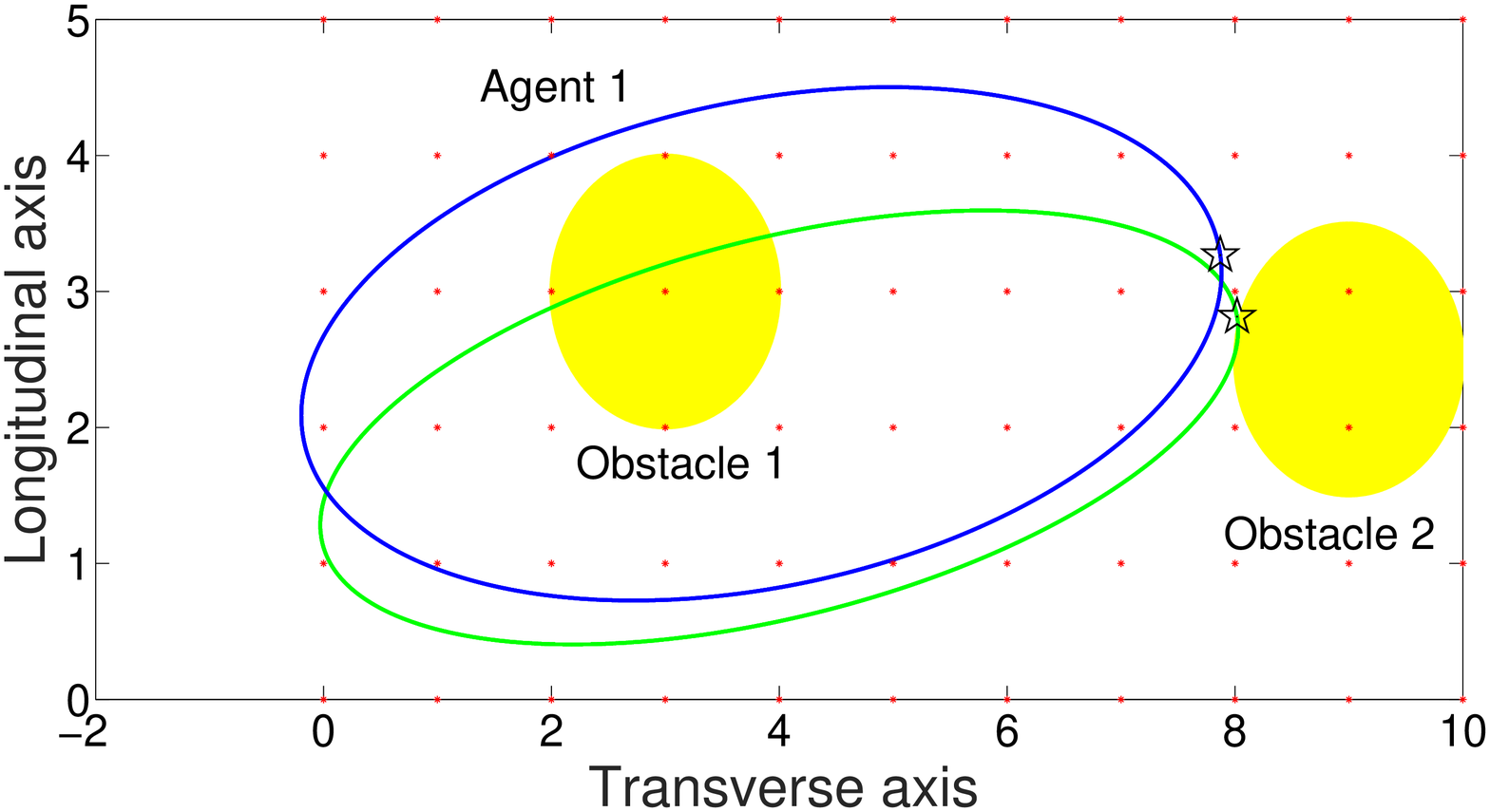}}
\subfigure[The evolution of performance metric J.]{\label{ZMJ_PersistentMonitoring_fig.62} 
\includegraphics[height=3.5cm,width=7cm]{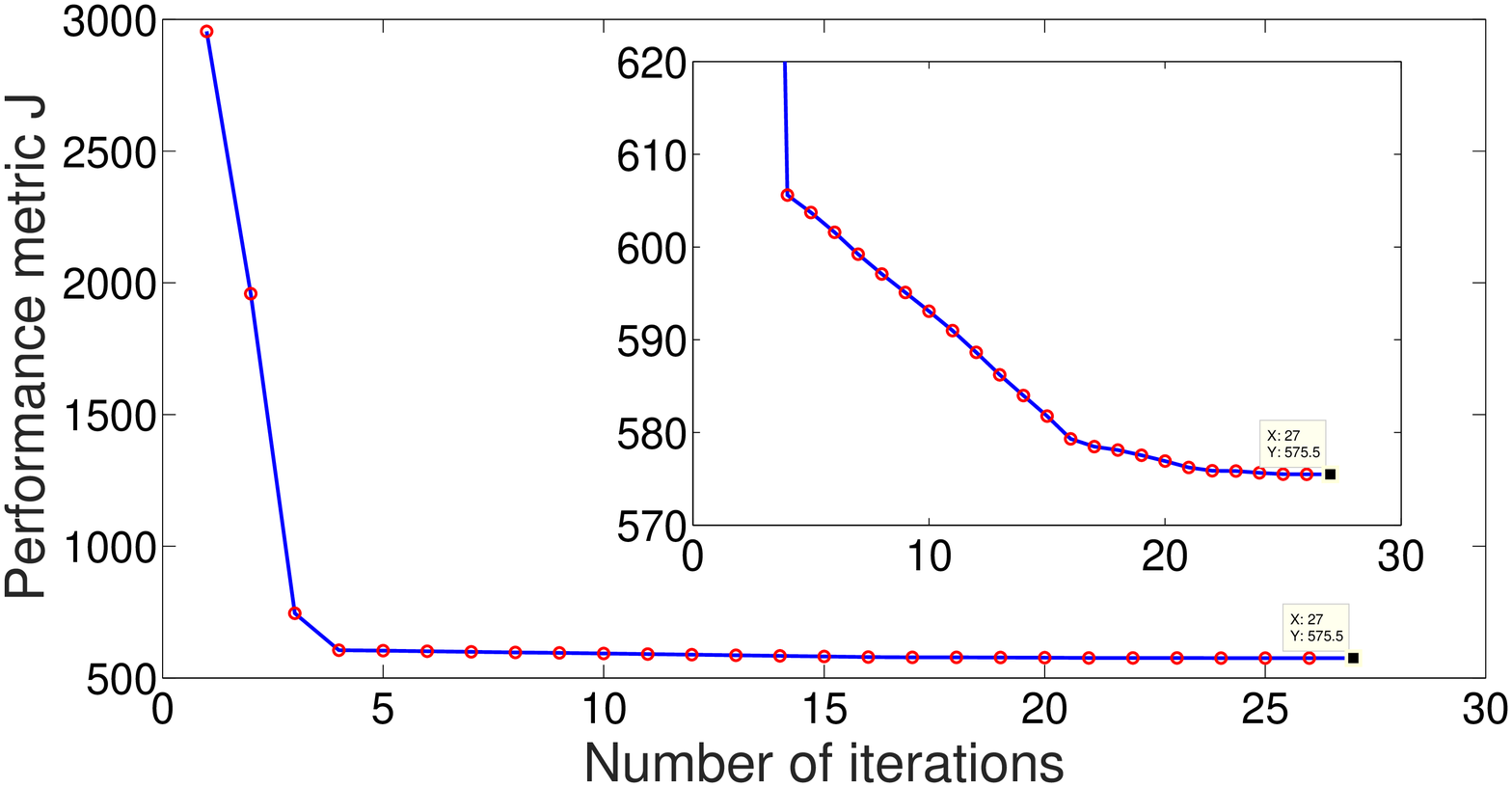}}
\subfigure[Distances between Agent 1 and Obstacles.]{\label{ZMJ_PersistentMonitoring_fig.63} 
\includegraphics[height=3.5cm,width=7cm]{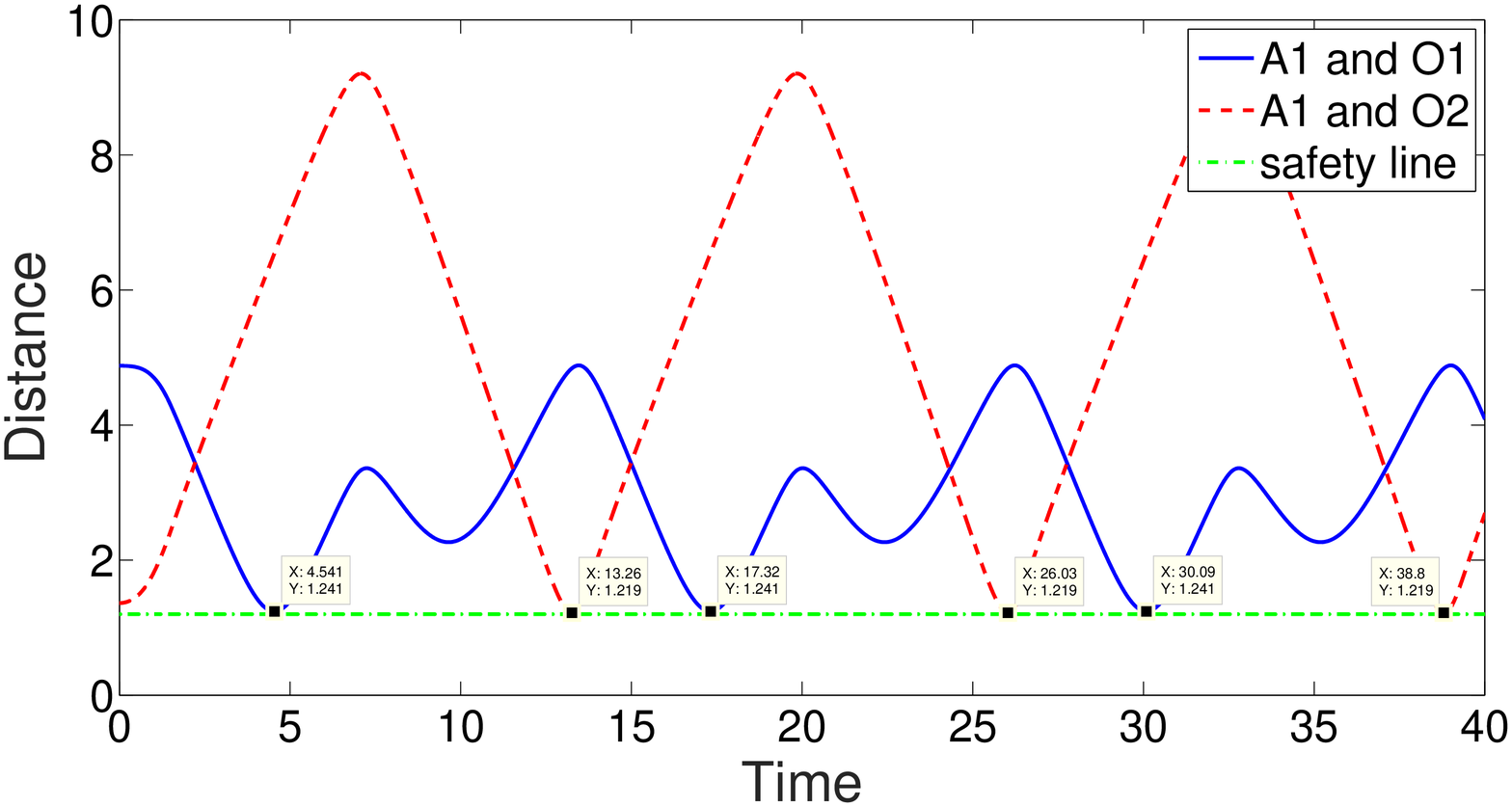}}
\caption{Elliptical trajectory: persistent monitoring task with obstacles using the IPA-based iteration algorithm for one agent (using the the sensing model in \cite{lin2015optimal}).} \label{ZMJ_PersistentMonitoring_fig.6}
\end{figure}

\textbf{\emph{Example 3.} Comparison with the existing sensing model.}

In this example, we use the existing sensing model in \cite{lin2015optimal}, which depends only on the relative distance between an agent and a target. Intuitively, we expect a better result than \textbf{\emph{Example 1}}, owing to the negative effect of moving agents in order to monitor targets. With the same initial conditions as in Fig. \ref{ZMJ_PersistentMonitoring_fig.4} of \textbf{\emph{Example 1}}, the results using the sensing model in \cite{lin2015optimal} are shown in Fig. \ref{ZMJ_PersistentMonitoring_fig.6}. Note that the final ellipse $[4.7291,2.2391,4.2087,1.6230,0.0321]$ in Fig. \ref{ZMJ_PersistentMonitoring_fig.61} is different from the final ellipse $[3.8791,2.4675,3.8994,1.8926,-0.0066]$ in Fig. \ref{ZMJ_PersistentMonitoring_fig.41}, which is due to the influence of the new sensing model on the final optimal trajectory. In addition, the value of the performance metric in Fig. \ref{ZMJ_PersistentMonitoring_fig.62} is less than that in Fig. \ref{ZMJ_PersistentMonitoring_fig.42}.
In Fig. \ref{ZMJ_PersistentMonitoring_fig.42}, the new sensing model is used which includes the effect of velocity on the quality of sensing. In Fig. \ref{ZMJ_PersistentMonitoring_fig.62} the original sensing model from \cite{lin2015optimal} is used. This ignores agent movements, therefore, it gives a more ``optimistic" result in the sense that the average uncertainty metric is lower than it ought to be with velocity taken into account. Actually, since the influence of the agent's velocity on sensing strength cannot be ignored, such result may be not achievable in practice. In addition, Fig. \ref{ZMJ_PersistentMonitoring_fig.63} can still verify that collisions are avoided through our method. Furthermore, we use the final ellipse in Fig. \ref{ZMJ_PersistentMonitoring_fig.61} to execute the monitoring task using the sensing model in this work, and the performance metric is $J = 662.6$, which is only slightly different from the result in Fig. \ref{ZMJ_PersistentMonitoring_fig.42}.

\textbf{\emph{Case B.} Two agents case.}

In this case, we carry out multiple simulation experiments and show the best results of two agents using the elliptical trajectory and the Fourier series trajectory to execute the monitoring task. The time horizon is set as $30s$. Suppose there are four targets $[5, 1], [5, 2], [5, 3], [5, 4]$ with larger weights $\sigma_{32} = \sigma_{33} = \sigma_{34} = \sigma_{35} = 2$ such that they require more attention than others. Moreover, obstacles in the mission space are covered by yellow circular areas whose centers are $[3, 3.8]$, $[8.5, 1.5]$, respectively and radii are both $1$.
\begin{figure}[!h]
\centering
\subfigure[Green ellipses: initial trajectories. Blue ellipses: final trajectories. Yellow circular areas: areas covering obstacles. Red $*$: target points. Black pentagrams: starting points.]{\label{ZMJ_PersistentMonitoring_fig.71} 
\includegraphics[height=3.5cm,width=7cm]{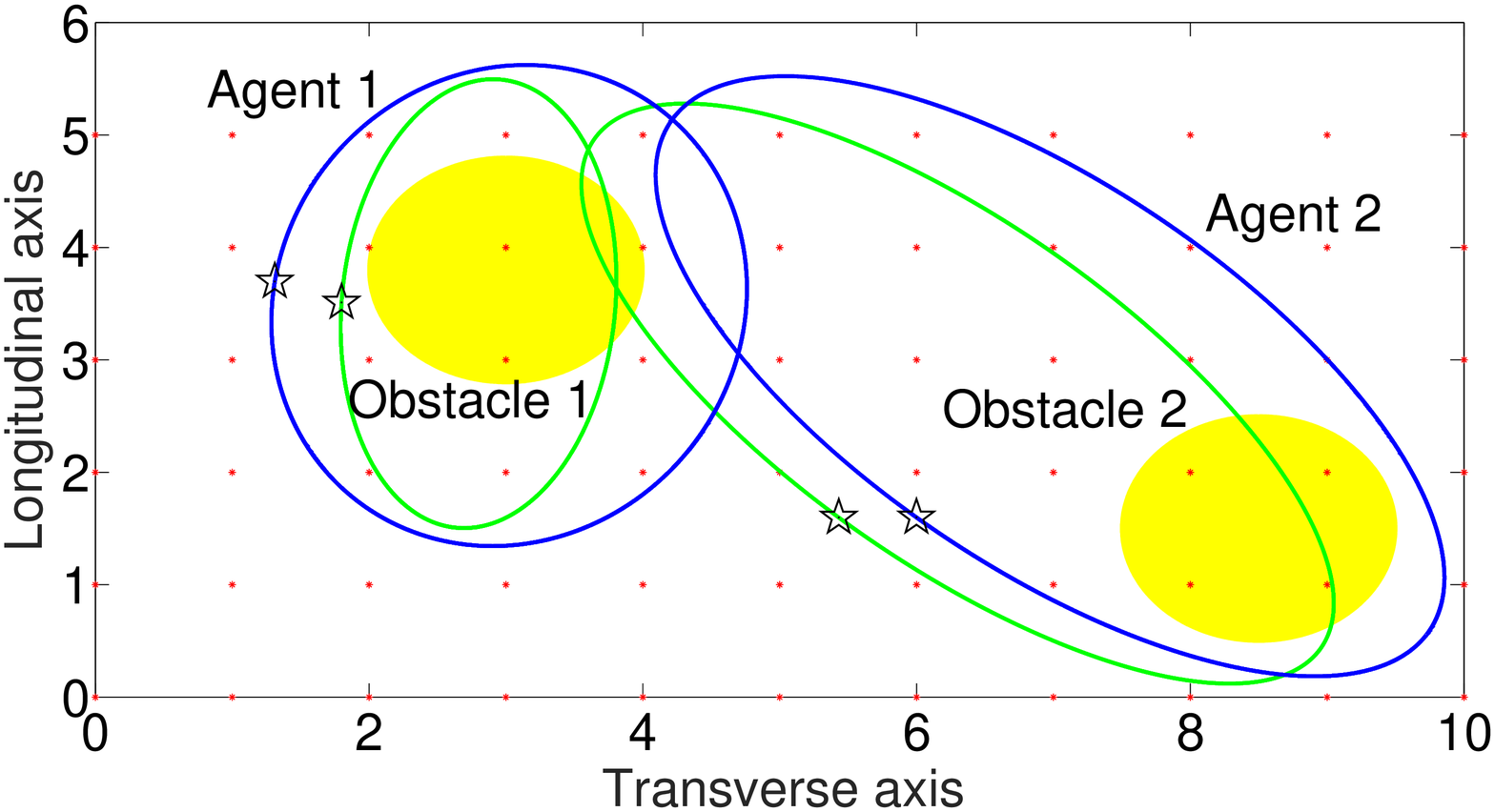}}
\subfigure[The evolution of performance metric J.]{\label{ZMJ_PersistentMonitoring_fig.72} 
\includegraphics[height=3.5cm,width=7cm]{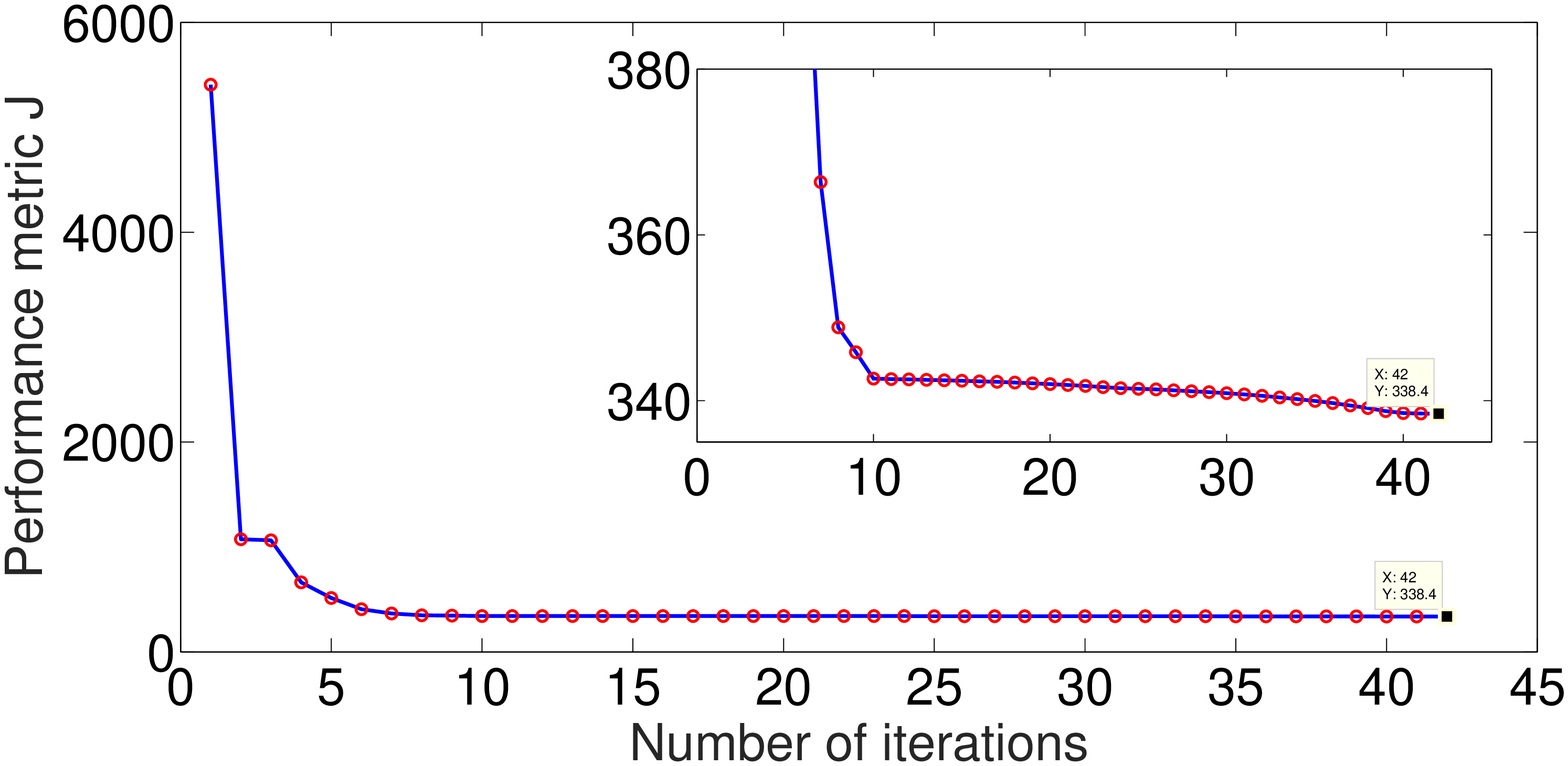}}
\subfigure[Distances between Agents and Obstacles.]{\label{ZMJ_PersistentMonitoring_fig.73} 
\includegraphics[height=3.5cm,width=7cm]{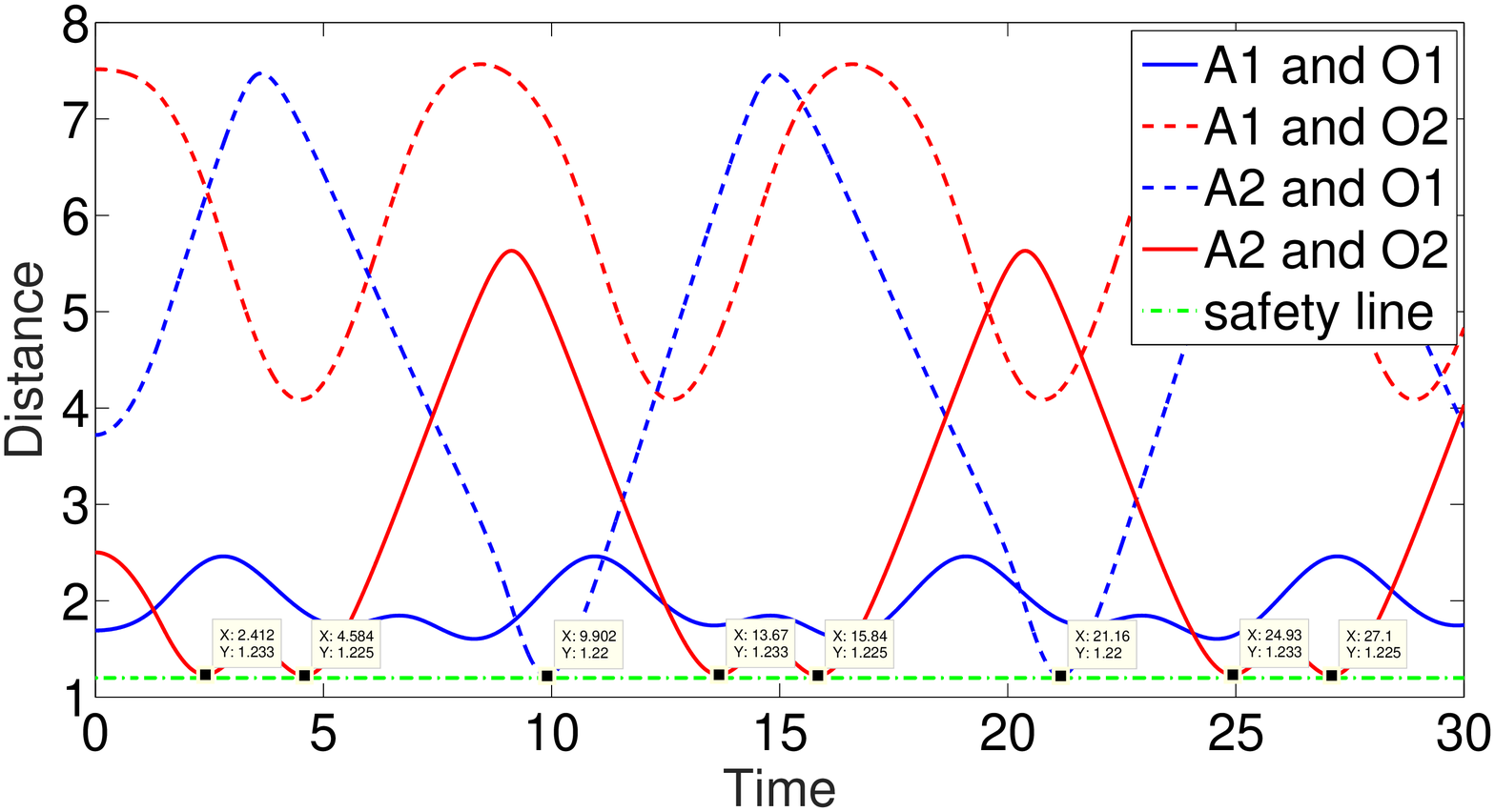}}
\subfigure[Distance between Agent 1 and Agent 2.]{\label{ZMJ_PersistentMonitoring_fig.74} 
\includegraphics[height=3.5cm,width=7cm]{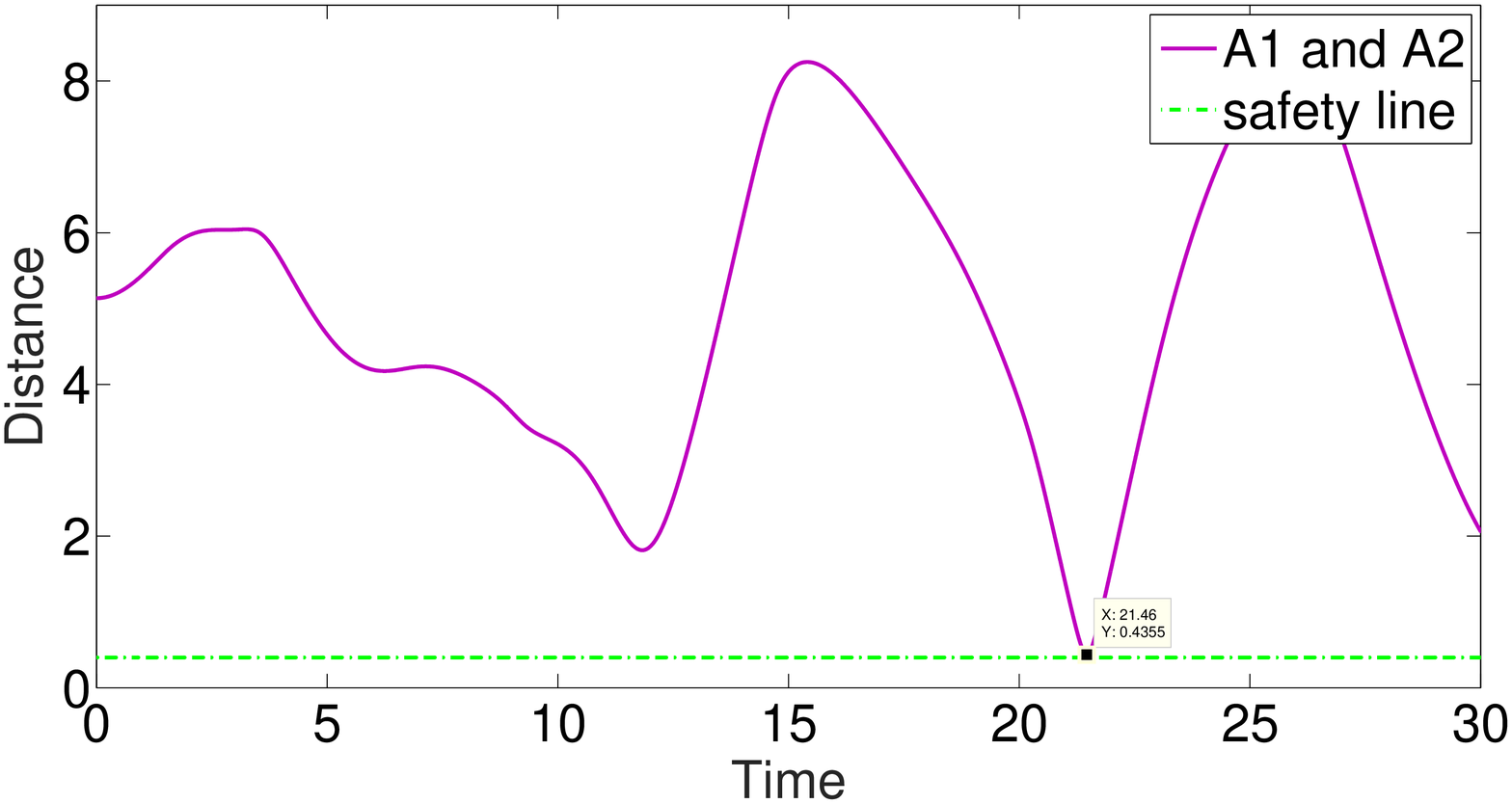}}
\caption{Elliptical trajectory: persistent monitoring task with obstacles using the IPA-based iteration algorithm for two agents.} \label{ZMJ_PersistentMonitoring_fig.7}
\end{figure}

Fig. \ref{ZMJ_PersistentMonitoring_fig.71} shows the optimal elliptical trajectories obtained by Algorithm \ref{ZMJ_PersistentMonitoring_algorithm1}. In Fig. \ref{ZMJ_PersistentMonitoring_fig.72}, the performance metric decreases as the number of iterations with the final value $J(\Theta^{42}) = 338.4$ and $|J(\Theta^{42}) - J(\Theta^{41})| < \varepsilon = 0.01$. Moreover, we show distances between agents and obstacles in Fig. \ref{ZMJ_PersistentMonitoring_fig.73}. It is obvious that distances between agents and obstacles are always greater than 1.2. Fig. \ref{ZMJ_PersistentMonitoring_fig.74} shows the distance between Agent 1 and Agent 2, which is always greater than 0.4. In other words, there are no potential collisions (agent to agent, agent to obstacle) in the final optimal trajectories.


Moreover, with the same setting as the elliptical example, simulation results of the Fourier trajectory example with two agents can be seen in Appendix C. And the final numerical results of the two trajectories are presented in TABLE \ref{ZMJ_PersistentMonitoring_table3}, and the final costs of collisions (agent to agent, agent to obstacle), i.e. $J_{2}^{\star}$ and $J_{3}^{\star}$, are zero, which indicates that we achieve collision-free monitoring.
\begin{table}[!h]
    \setlength{\abovecaptionskip}{-0.2cm}
    \setlength{\belowcaptionskip}{-0cm}
    \renewcommand\arraystretch{1.3}
    \caption{Result Comparison For Case B}\label{ZMJ_PersistentMonitoring_table3}
    \begin{center}
        \begin{tabular}{|c|c|c|c|}
            \hline
            $Method$ & $J^{\star}$ & $J_{2}^{\star}$ & $J_{3}^{\star}$ \\
            \hline
            $Elliptical$ & $338.4$ & $0$ & $0$ \\
            \hline
            $Fourier$ & $305.9$ & $0$ & $0$ \\
            \hline
       \end{tabular}
    \end{center}
\end{table}

\section{Conclusion}\label{ZMJ_PersistentMonitoring_section_6}
We have studied a 2D persistent monitoring problem using second-order agents. We have proposed a new and more practical sensing model which takes the effects of both distances between agents and targets and agents' velocities into account. In particular, we have  considered no-collision constraints (both agent to agent and agent to obstacle). Through  parameterizing agent trajectories and utilizing IPA-based gradient algorithm, our method provides collision-free optimal (locally) trajectories.

In addition, decentralizing the controller is particularly challenging in persistent monitoring. We have been working on this problem and report results in \cite{zhou2017decentralized}. Further, in \cite{zhou2019the} we show that it is possible to achieve an ``almost decentralized" solution based exclusively on using local information. We also quantify the ``price of decentralization" in the sense of losing some performance. Thus, imposing a decentralized solution to the persistent monitoring problem is of course possible but not without some (quantifiable) performance degradation. Future research will focus on escaping the local optimal solution, developing decentralized solutions, developing nonconservative ways of dealing with obstacles and exploiting on-line strategies of dealing with unexpected events for the 2D persistent monitoring problems with potential collisions.


%

%

%
\appendices
\section{Proof of Proposition 1}\label{ZMJ_PersistentMonitoring_appendixA}
The proof will be divided into two parts based on the two types of events.

(i) Event $\xi_{i}^{0}$: This is an endogenous event which causes a transition of \eqref{ZMJ_PersistentMonitoring_Uncertainty} from $R_{i}(t) > 0, t < \rho_{k}$ to $R_{i}(t) = 0, t \geq \rho_{k}$. Set $g_{k}(\Theta,\bm{\chi}) = R_{i}(\rho_{k}) = 0$. Then from (\ref{ZMJ_PersistentMonitoring_IPA3}) and (\ref{ZMJ_PersistentMonitoring_Uncertainty}), we can get
\begin{equation}
\begin{aligned}
\rho_{k} '
&= -[\frac{\partial g_{k}}{\partial R_{i}}f_{k}(\rho_{k}^{-})]^{-1}[\frac{\partial g_{k}}{\partial \bm{\Theta}} + \frac{\partial g_{k}}{\partial R_{i}}R_{i} '(\rho_{k}^{-})]\\
&= -\frac{R_{i} '(\rho_{k}^{-})}{A_{i} - BP_{i}(\rho_{k}^{-})}
\end{aligned}
\end{equation}
Following (\ref{ZMJ_PersistentMonitoring_IPA2}) we have
\begin{equation}
\begin{aligned}
\nabla R_{i}(\rho_{k}^{+}) = \nabla R_{i}(\rho_{k}^{-}) -\frac{[A_{i}-BP_{i}(\rho_{k}^{-})]R_{i} '(\rho_{k}^{-})}{A_{i} - BP_{i}(\rho_{k}^{-})} = 0\\
\end{aligned}
\end{equation}
Therefore, $\nabla R_{i}(\rho_{k}^{+})$ is reset to 0 if event $\xi_{i}^{0}$ happens.

(ii) Event $\xi_{i}^{+}$: This event causes a transition of \eqref{ZMJ_PersistentMonitoring_Uncertainty} from $R_{i}(t) = 0, t \leq \rho_{k}$ to $R_{i}(t) > 0, t > \rho_{k}$. In this case, at $\rho_{k}$ we have $A_{i}(\rho_{k}) = BP_{i}(\rho_{k})$. Obviously, $R_{i}(t)$ is continuous and $f_{k-1}(\rho_{k}^{-}) = f_{k}(\rho_{k}^{+})$ in (\ref{ZMJ_PersistentMonitoring_IPA2}). Therefore, we get
\begin{equation}
\begin{aligned}
\nabla R_{i}(\rho_{k}^{+}) = \nabla R_{i}(\rho_{k}^{-})
\end{aligned}
\end{equation}
under event $\xi_{i}^{+}$.
$\hfill\blacksquare$

\section{Derivations of Fourier Trajectories}\label{ZMJ_PersistentMonitoring_appendixB}
In order to further enrich the work of trajectory optimization, we also select the Fourier series trajectories \cite{khazaeni2018event,zahn1972fourier} to solve the persistent monitoring problem. And then there is the derivation results about Fourier series trajectories.

The agent $n$'s trajectory thus can be described as follows,
\begin{equation}\label{ZMJ_PersistentMonitoring_s_Fourier}
\begin{aligned}
s_{n}^{x}(t) &= a_{n,0} + \displaystyle{\sum_{\gamma = 1}^{\Gamma_{n}^{x}}}a_{n,\gamma}\sin(2\pi\gamma f_{n}^{x}\rho_{n}(t) + \phi_{n,\gamma}^{x})\\
s_{n}^{y}(t) &= b_{n,0} + \displaystyle{\sum_{\gamma = 1}^{\Gamma_{n}^{y}}}b_{n,\gamma}\sin(2\pi\gamma f_{n}^{y}\rho_{n}(t) + \phi_{n,\gamma}^{y})
\end{aligned}
\end{equation}
where $f_{n}^{x}$ and$f_{n}^{y}$ are the base frequencies, $a_{n,0}$ and $b_{n,0}$ are the zero frequency components, $a_{n,\gamma}$ and $b_{n,\gamma}$ are the amplitudes for the sinusoid functions, $\phi_{n,\gamma}^{x}$ and $\phi_{n,\gamma}^{y}$ are the phase differences. And $\rho_{n}(t) \in [0,2\pi)$ is the agent position on the trajectory. Then by taking (\ref{ZMJ_PersistentMonitoring_s_Fourier}), we can get the velocity information of the agent as follows,
\begin{equation}\label{ZMJ_PersistentMonitoring_v_Fourier}
\begin{aligned}
\dot{s}_{n}^{x}(t) &= 2\pi f_{n}^{x}\dot{\rho}_{n}(t)\displaystyle{\sum_{\gamma = 1}^{\Gamma_{n}^{x}}}\gamma a_{n,\gamma}\cos(2\pi\gamma f_{n}^{x}\rho_{n}(t) + \phi_{n,\gamma}^{x})\\
\dot{s}_{n}^{y}(t) &= 2\pi f_{n}^{y}\dot{\rho}_{n}(t)\displaystyle{\sum_{\gamma = 1}^{\Gamma_{n}^{y}}}\gamma b_{n,\gamma}\cos(2\pi\gamma f_{n}^{y}\rho_{n}(t) + \phi_{n,\gamma}^{y})
\end{aligned}
\end{equation}
Further taking the derivative of (\ref{ZMJ_PersistentMonitoring_v_Fourier}), we can get the acceleration information of the agent as follows,
\begin{equation}\label{ZMJ_PersistentMonitoring_u_Fourier}
\begin{aligned}
\ddot{s}_{n}^{x}(t) &= 2\pi f_{n}^{x}\ddot{\rho}_{n}(t)\displaystyle{\sum_{\gamma = 1}^{\Gamma_{n}^{x}}}\gamma a_{n,\gamma}\cos(2\pi\gamma f_{n}^{x}\rho_{n}(t) + \phi_{n,\gamma}^{x})\\
&- (2\pi f_{n}^{x}\dot{\rho}_{n}(t))^2\displaystyle{\sum_{\gamma = 1}^{\Gamma_{n}^{x}}}\gamma^{2} a_{n,\gamma}\sin(2\pi\gamma f_{n}^{x}\rho_{n}(t) + \phi_{n,\gamma}^{x})\\
\ddot{s}_{n}^{y}(t) &= 2\pi f_{n}^{y}\ddot{\rho}_{n}(t)\displaystyle{\sum_{\gamma = 1}^{\Gamma_{n}^{y}}}\gamma b_{n,\gamma}\cos(2\pi\gamma f_{n}^{y}\rho_{n}(t) + \phi_{n,\gamma}^{y})\\
&- (2\pi f_{n}^{y}\dot{\rho}_{n}(t))^2\displaystyle{\sum_{\gamma = 1}^{\Gamma_{n}^{y}}}\gamma^{2} b_{n,\gamma}\sin(2\pi\gamma f_{n}^{y}\rho_{n}(t) + \phi_{n,\gamma}^{y})
\end{aligned}
\end{equation}

For the Fourier series trajectories, the parameter of agent $n$ is $\Theta_{n} = [f_{n}^{x},f_{n}^{y},a_{n,0},...,a_{n,\Gamma_{n}^{x}},b_{n,0},...,b_{n,\Gamma_{n}^{y}},\phi_{n,1}^{x},...,$
$\phi_{n,\Gamma_{n}^{x}}^{x},\phi_{n,1}^{y},...,\phi_{n,\Gamma_{n}^{y}}^{y}]$. Due to $\frac{f_{n}^{x}}{f_{n}^{y}}$ determining the shape of the trajectory, therefore, we will keep $f_{n}^{y}$ constant and adjust $f_{n}^{x}$. Then the gradients about these parameters can be listed as follows,
\begin{equation}
\begin{aligned}
\frac{\partial s_{n}^{x}}{\partial f_{n}^{x}} = 2\pi \rho_{n}(t)\displaystyle{\sum_{\gamma = 1}^{\Gamma_{n}^{x}}}\gamma a_{n,\gamma}\cos(2\pi\gamma f_{n}^{x}\rho_{n}(t) + \phi_{n,\gamma}^{x})
\end{aligned}
\end{equation}
\begin{equation}
\begin{aligned}
\frac{\partial s_{n}^{x}}{\partial a_{n,0}} = 1,\qquad \frac{\partial s_{n}^{x}}{\partial b_{n,0}} = 0
\end{aligned}
\end{equation}
\begin{equation}
\begin{aligned}
\frac{\partial s_{n}^{x}}{\partial a_{n,\gamma}} = \sin(2\pi\gamma f_{n}^{x}\rho_{n}(t) + \phi_{n,\gamma}^{x}),\qquad \frac{\partial s_{n}^{x}}{\partial b_{n,\gamma}} = 0
\end{aligned}
\end{equation}
\begin{equation}
\begin{aligned}
\frac{\partial s_{n}^{x}}{\partial \phi_{n,\gamma}^{x}} = a_{n,\gamma}\cos(2\pi\gamma f_{n}^{x}\rho_{n}(t) + \phi_{n,\gamma}^{x}),\quad \frac{\partial s_{n}^{x}}{\partial \phi_{n,\gamma}^{y}} = 0
\end{aligned}
\end{equation}
\begin{equation}
\begin{aligned}
\frac{\partial s_{n}^{y}}{\partial f_{n}^{x}} = 0,\qquad \frac{\partial s_{n}^{y}}{\partial a_{n,0}} = 0,\qquad \frac{\partial s_{n}^{y}}{\partial b_{n,0}} = 1
\end{aligned}
\end{equation}
\begin{equation}
\begin{aligned}
\frac{\partial s_{n}^{y}}{\partial a_{n,\gamma}} = 0,\qquad \frac{\partial s_{n}^{y}}{\partial b_{n,\gamma}} = \sin(2\pi\gamma f_{n}^{y}\rho_{n}(t) + \phi_{n,\gamma}^{y})
\end{aligned}
\end{equation}
\begin{equation}
\begin{aligned}
\frac{\partial s_{n}^{y}}{\partial \phi_{n,\gamma}^{x}} = 0,\quad \frac{\partial s_{n}^{y}}{\partial \phi_{n,\gamma}^{y}} = b_{n,\gamma}\cos(2\pi\gamma f_{n}^{y}\rho_{n}(t) + \phi_{n,\gamma}^{y})
\end{aligned}
\end{equation}
\begin{equation}
\begin{aligned}
\frac{\partial v_{n}^{x}}{\partial f_{n}^{x}} = 2\pi \dot{\rho}_{n}(t)\displaystyle{\sum_{\gamma = 1}^{\Gamma_{n}^{x}}}\gamma a_{n,\gamma}\cos(2\pi\gamma f_{n}^{x}\rho_{n}(t) + \phi_{n,\gamma}^{x})
\end{aligned}
\end{equation}
\begin{equation}
\begin{aligned}
\frac{\partial v_{n}^{x}}{\partial a_{n,0}} = 0,\qquad \frac{\partial v_{n}^{x}}{\partial b_{n,0}} = 0
\end{aligned}
\end{equation}
\begin{equation}
\begin{aligned}
\frac{\partial v_{n}^{x}}{\partial a_{n,\gamma}} = 2\pi f_{n}^{x}\dot{\rho}_{n}(t)\gamma \cos(2\pi\gamma f_{n}^{x}\rho_{n}(t) + \phi_{n,\gamma}^{x}), \frac{\partial v_{n}^{x}}{\partial b_{n,\gamma}} = 0
\end{aligned}
\end{equation}
\begin{equation}
\begin{aligned}
\frac{\partial v_{n}^{x}}{\partial \phi_{n,\gamma}^{x}} &= -2\pi f_{n}^{x}\dot{\rho}_{n}(t)\gamma a_{n,\gamma}\sin(2\pi\gamma f_{n}^{x}\rho_{n}(t) + \phi_{n,\gamma}^{x}),\\
\frac{\partial v_{n}^{x}}{\partial \phi_{n,\gamma}^{y}} &= 0
\end{aligned}
\end{equation}
\begin{equation}
\begin{aligned}
\frac{\partial v_{n}^{y}}{\partial f_{n}^{x}} = 0,\qquad \frac{\partial v_{n}^{y}}{\partial a_{n,0}} = 0,\qquad \frac{\partial v_{n}^{y}}{\partial b_{n,0}} = 0
\end{aligned}
\end{equation}
\begin{equation}
\begin{aligned}
\frac{\partial v_{n}^{y}}{\partial a_{n,\gamma}} = 0, \frac{\partial v_{n}^{y}}{\partial b_{n,\gamma}} = 2\pi f_{n}^{y}\dot{\rho}_{n}(t)\gamma \cos(2\pi\gamma f_{n}^{y}\rho_{n}(t) + \phi_{n,\gamma}^{y})
\end{aligned}
\end{equation}
\begin{equation}
\begin{aligned}
\frac{\partial v_{n}^{y}}{\partial \phi_{n,\gamma}^{x}} &= 0,\\
\frac{\partial v_{n}^{y}}{\partial \phi_{n,\gamma}^{y}} &= -2\pi f_{n}^{y}\dot{\rho}_{n}(t)\gamma b_{n,\gamma}\sin(2\pi\gamma f_{n}^{y}\rho_{n}(t) + \phi_{n,\gamma}^{y})
\end{aligned}
\end{equation}

\section{Simulations of Fourier Trajectories}\label{ZMJ_PersistentMonitoring_appendixC}
\begin{figure}[!t]
\centering
\subfigure[Green Fourier trajectory: initial trajectory. Blue Fourier trajectory: final trajectory. Yellow circular areas: areas covering obstacles. Red $*$: target points. Black pentagrams: starting points.]{\label{ZMJ_PersistentMonitoring_fig.F1S} 
\includegraphics[height=3.5cm,width=7cm]{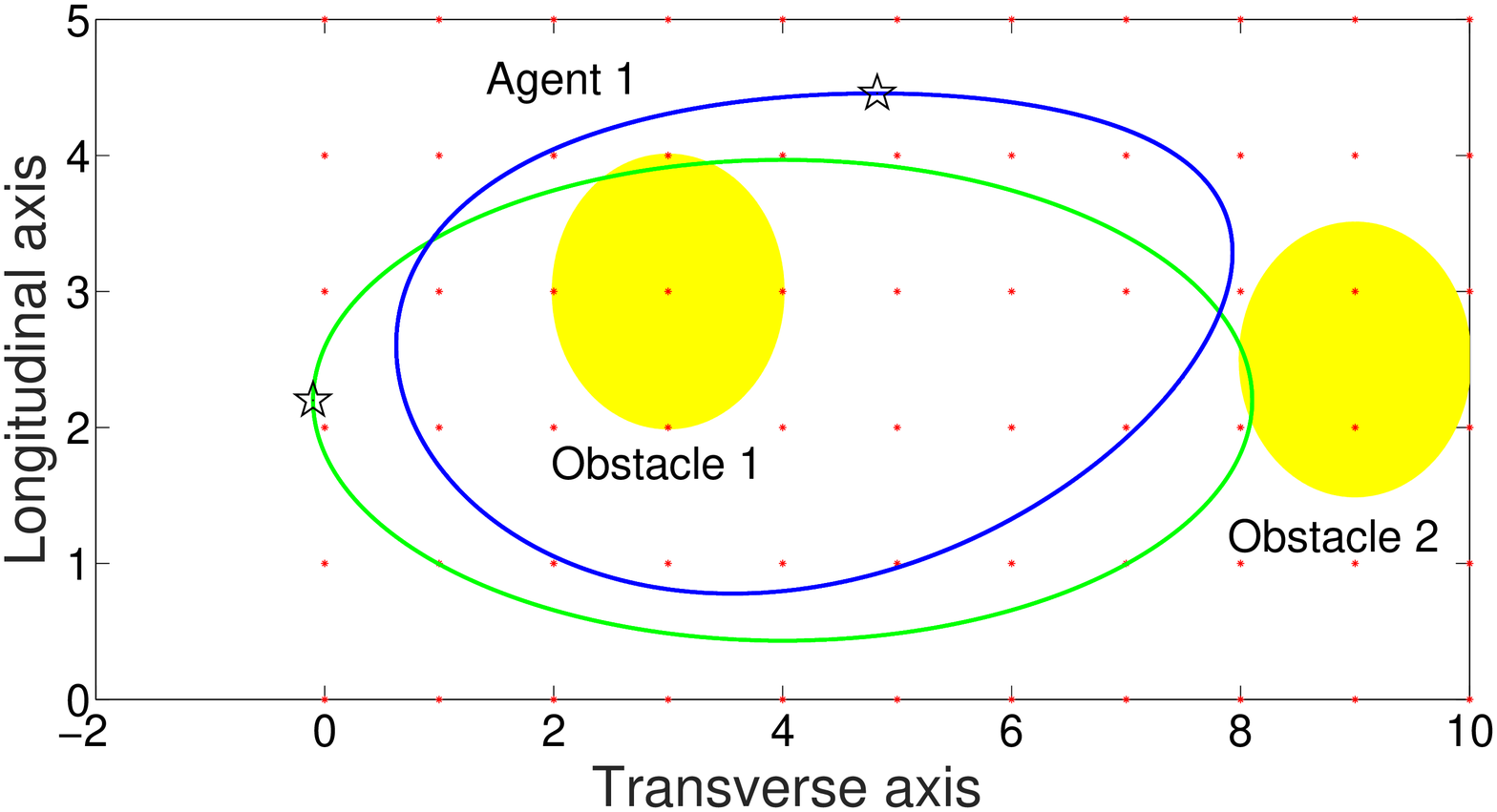}}
\subfigure[The evolution of performance metric J.]{\label{ZMJ_PersistentMonitoring_fig.F1J} 
\includegraphics[height=3.5cm,width=7cm]{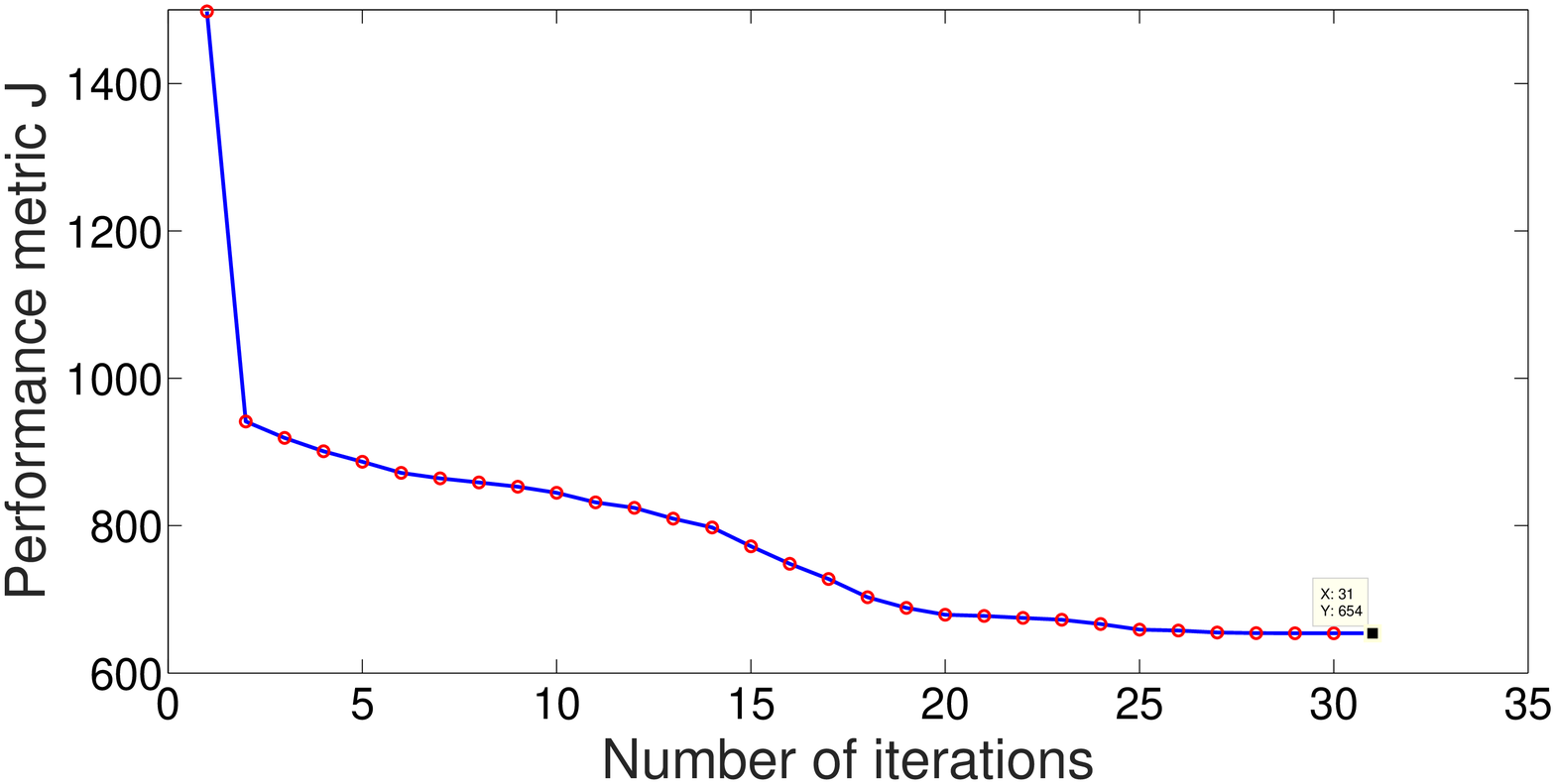}}
\subfigure[Distances between Agent 1 and Obstacles.]{\label{ZMJ_PersistentMonitoring_fig.F1dis} 
\includegraphics[height=3.5cm,width=7cm]{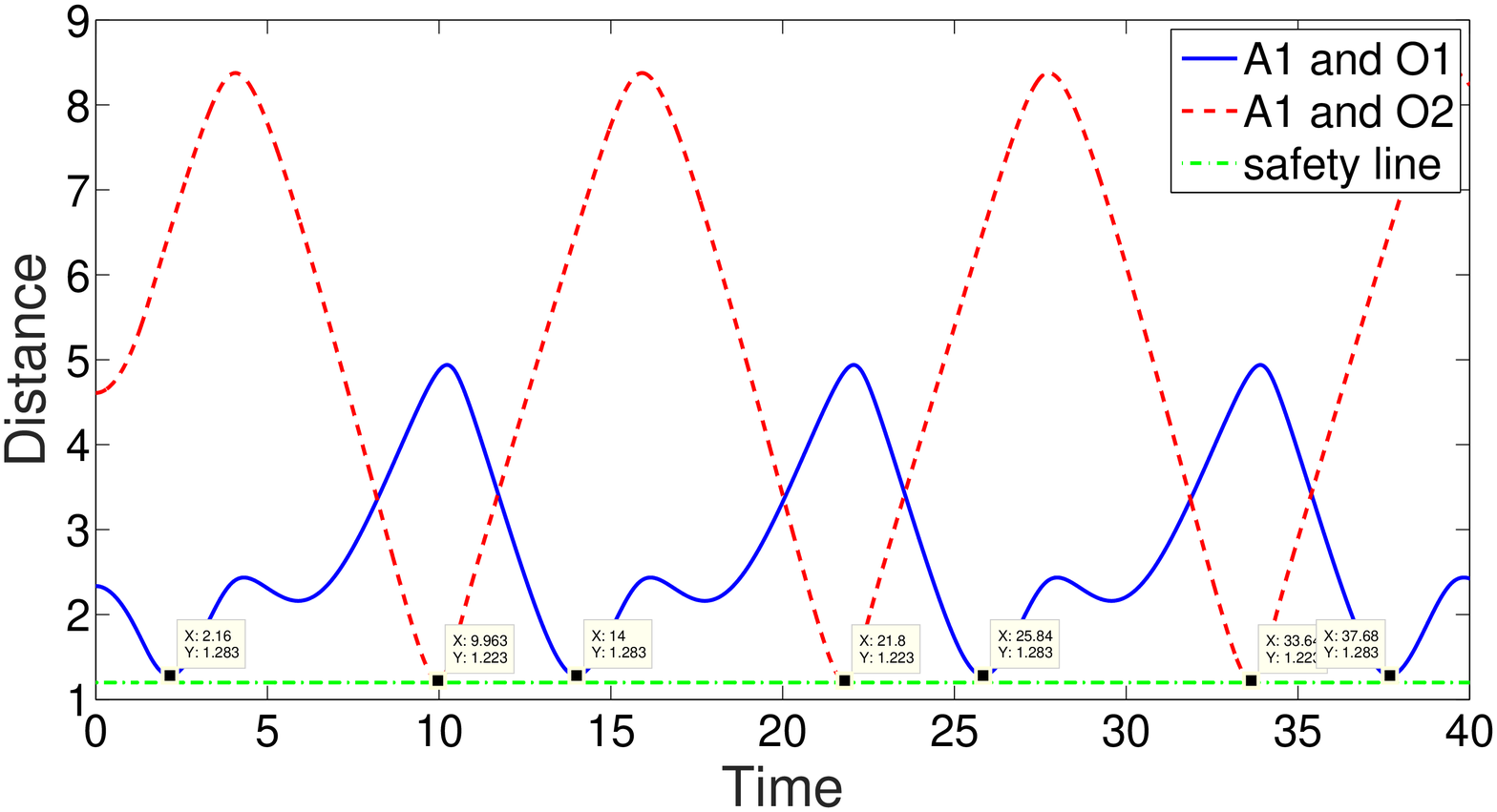}}
\caption{Fourier series trajectory: persistent monitoring task with obstacles using the IPA-based iteration algorithm for one agent.} \label{ZMJ_PersistentMonitoring_fig.F11}
\end{figure}

\begin{figure}[!t]
\centering
\subfigure[Green Fourier trajectories: initial trajectories. Blue Fourier trajectories: final trajectories. Yellow circular areas: areas covering obstacles. Red $*$: target points. Black pentagrams: starting points.]{\label{ZMJ_PersistentMonitoring_fig.F1} 
\includegraphics[height=3.5cm,width=7cm]{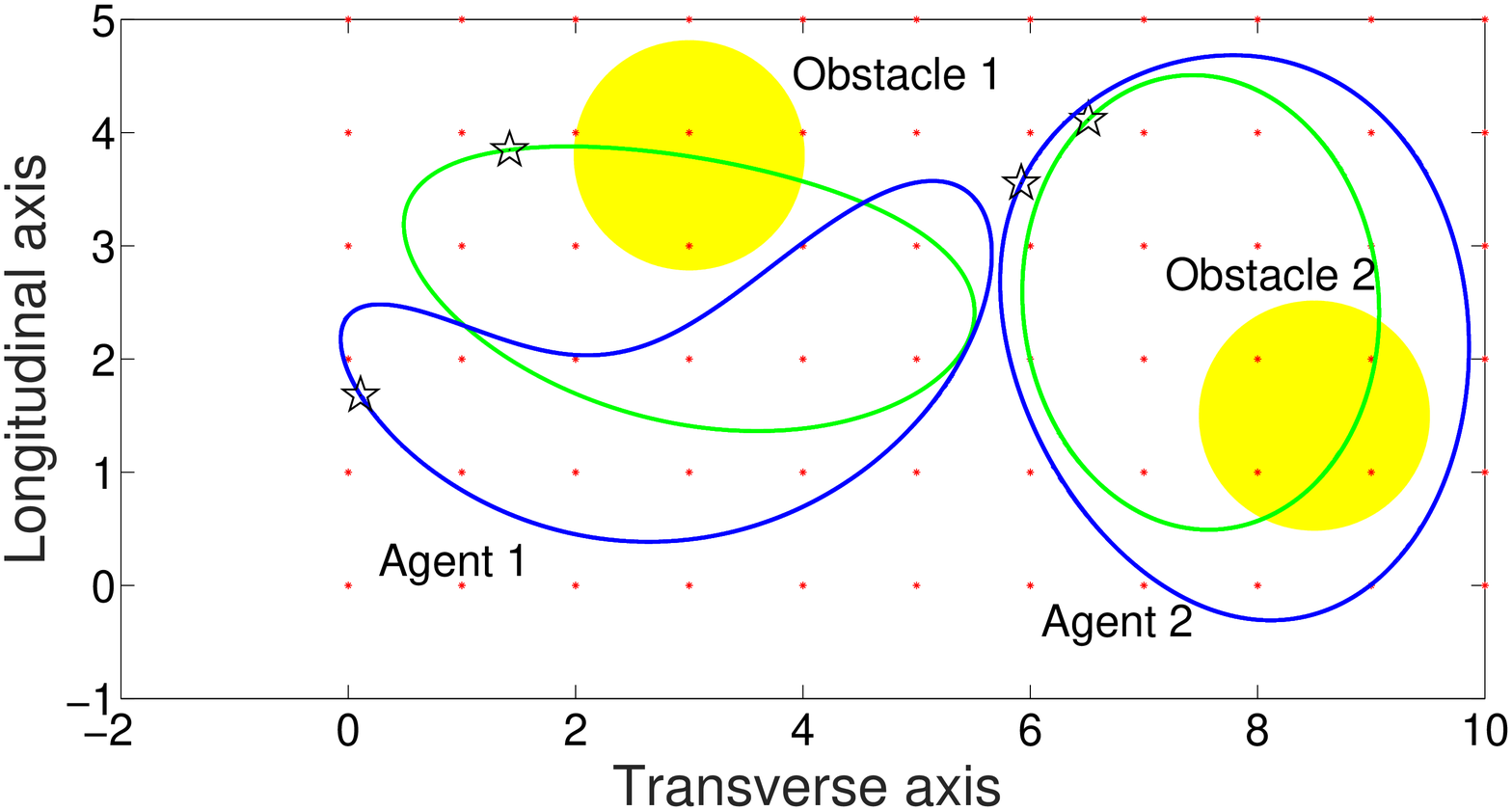}}
\subfigure[The evolution of performance metric J.]{\label{ZMJ_PersistentMonitoring_fig.F2} 
\includegraphics[height=3.5cm,width=7cm]{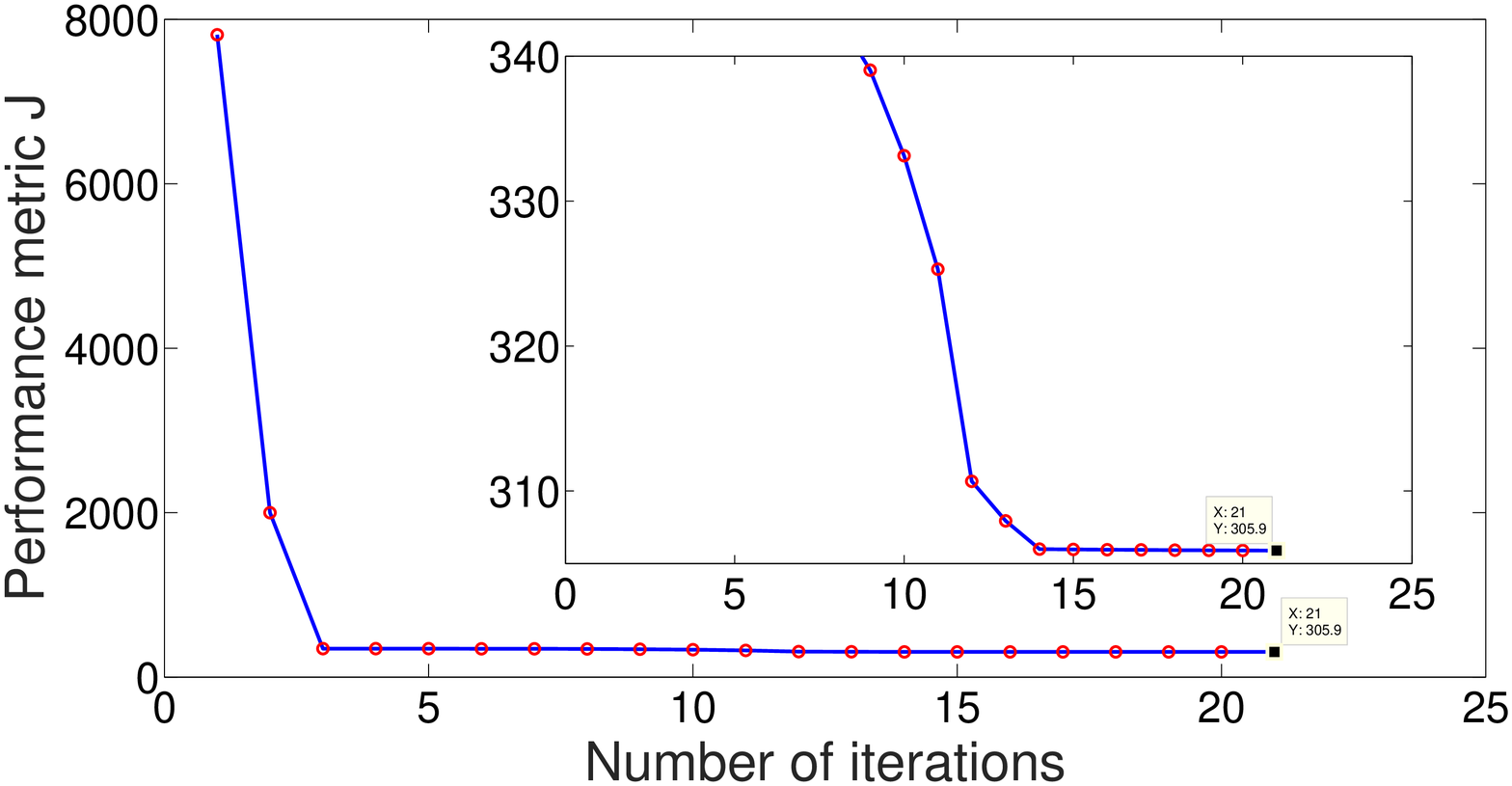}}
\subfigure[Distances between Agents and Obstacles.]{\label{ZMJ_PersistentMonitoring_fig.F3} 
\includegraphics[height=3.5cm,width=7cm]{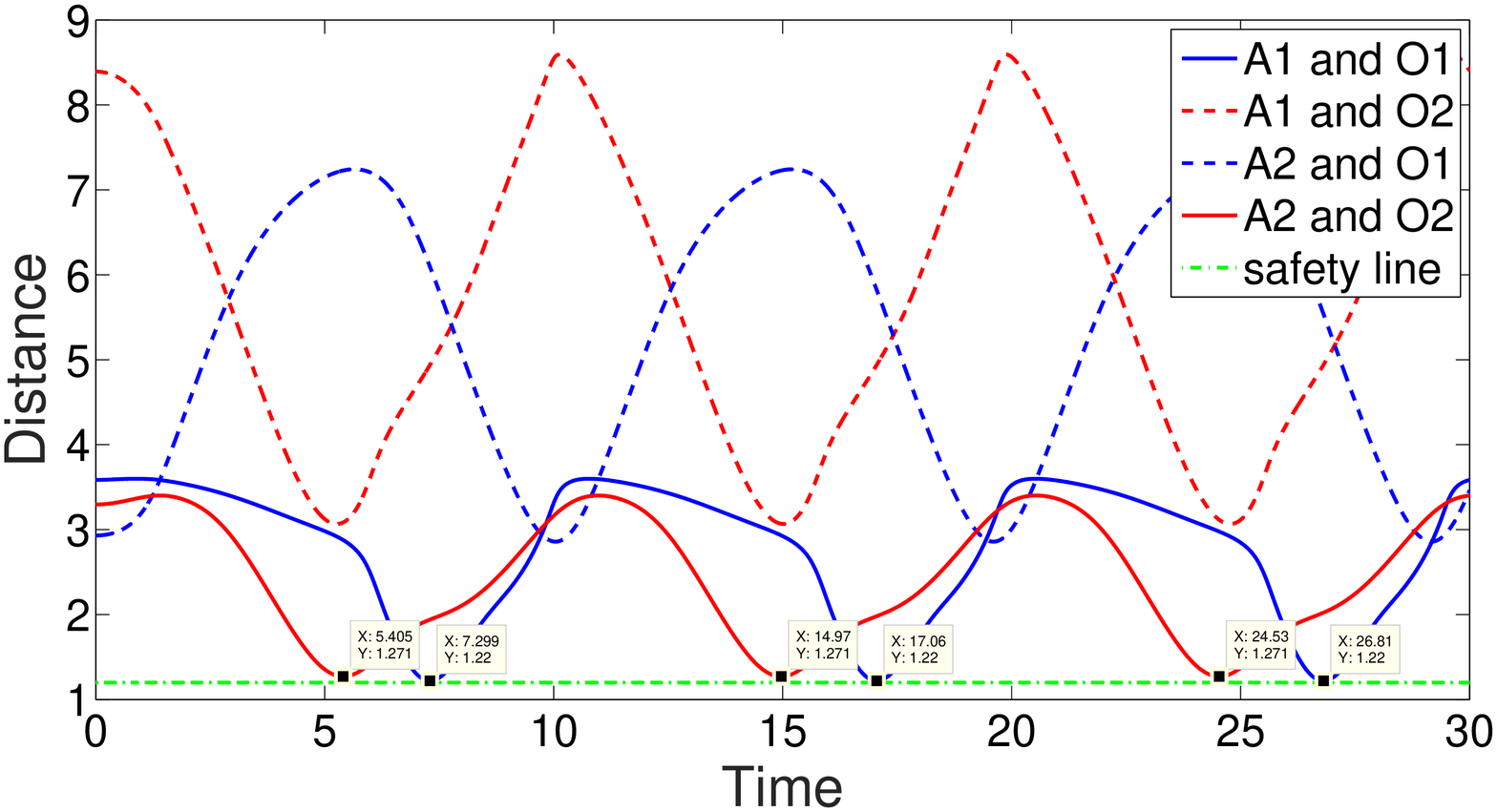}}
\subfigure[Distance between Agent 1 and Agent 2.]{\label{ZMJ_PersistentMonitoring_fig.F4} 
\includegraphics[height=3.5cm,width=7cm]{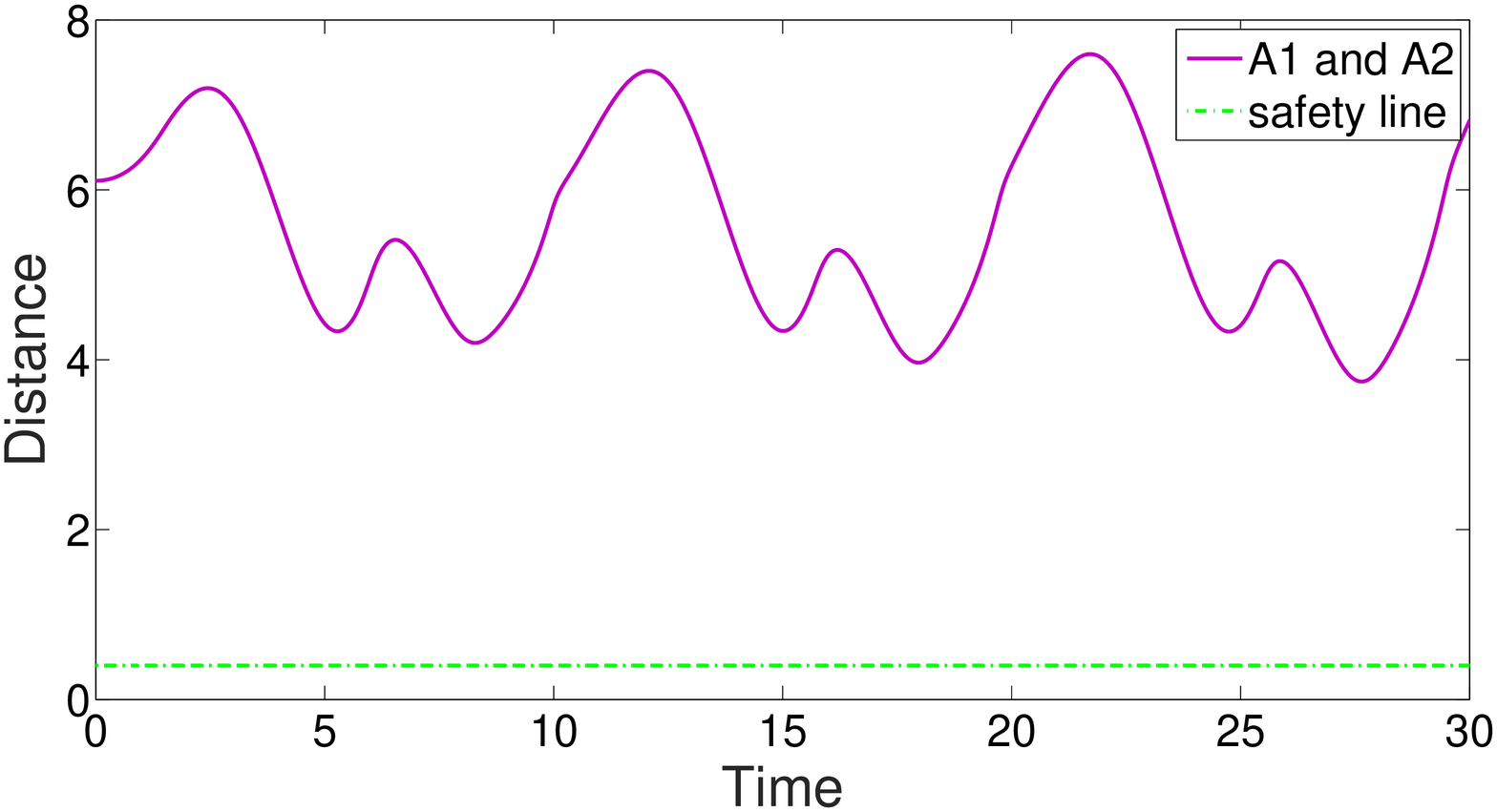}}
\caption{Fourier series trajectory: persistent monitoring task with obstacles using the IPA-based iteration algorithm for two agents.} \label{ZMJ_PersistentMonitoring_fig.F}
\end{figure}

With the same setting as the elliptical trajectory, the simulation results of the Fourier trajectory are shown in Fig. \ref{ZMJ_PersistentMonitoring_fig.F11} and Fig. \ref{ZMJ_PersistentMonitoring_fig.F}, respectively.

In Fig. \ref{ZMJ_PersistentMonitoring_fig.F11}, there is a persistent monitoring task with obstacles executed by one agent moving on a Fourier series trajectory and we select the best from multiple Fourier results. Because of the complexity and variability of this trajectory, we used $\Gamma_{1}^{x}=\Gamma_{1}^{y}=2$ for simplicity. From Fig. \ref{ZMJ_PersistentMonitoring_fig.F1S}, there is an interesting difference between the initial and final trajectory, which reflects the unpredictability of this trajectory. In Fig. \ref{ZMJ_PersistentMonitoring_fig.F1J}, the performance metric decreases as the iteration progresses and ultimately converges, thus the effectiveness of Algorithm \ref{ZMJ_PersistentMonitoring_algorithm1} applied to a Fourier series trajectory is verified. The final performance metric $J(\Theta^{31}) = 654$ and $|J(\Theta^{31}) - J(\Theta^{30})| < \varepsilon = 0.01$. From the distances between Agent 1 and obstacles in Fig. \ref{ZMJ_PersistentMonitoring_fig.F1dis}, we can observe that Agent 1 safely completed the monitoring task.

In Fig. \ref{ZMJ_PersistentMonitoring_fig.F}, there is a persistent monitoring task with obstacles executed by two agents moving on Fourier series trajectories. We set $\Gamma_{1}^{x}=\Gamma_{1}^{y}=\Gamma_{2}^{x}=\Gamma_{2}^{y}=2$ for simplicity. Fig. \ref{ZMJ_PersistentMonitoring_fig.F1} shows the interesting changes between the initial and final trajectories, which is due to agents tending to avoid obstacles. In Fig. \ref{ZMJ_PersistentMonitoring_fig.F2}, the performance metric decreases as the iteration progresses with $J(\Theta^{21}) = 305.9$ and $|J(\Theta^{21}) - J(\Theta^{20})| < \varepsilon = 0.01$. From Fig. \ref{ZMJ_PersistentMonitoring_fig.F3}, we can observe that the two agents both can safely execute the monitoring task. And in Fig. \ref{ZMJ_PersistentMonitoring_fig.F4}, the distance between Agent 1 and Agent 2 is always greater than 0.4, i.e. this is a safe monitoring task.

\ifCLASSOPTIONcaptionsoff
  \newpage
\fi

\bibliographystyle{ieeetr}        
\bibliography{references}

\end{document}